\documentclass{article}

\usepackage{arxiv}

\usepackage[utf8]{inputenc} 
\usepackage[T1]{fontenc}    
\usepackage{hyperref}       
\usepackage{url}            
\usepackage{booktabs}       
\usepackage{amsfonts}       
\usepackage{nicefrac}       
\usepackage{microtype}      
\usepackage{lipsum}
\usepackage{graphicx}

\usepackage{natbib}
\usepackage{subfig}
\usepackage{acro}
\usepackage{geometry}
\usepackage{comment}
\usepackage{xspace}
\usepackage{amsmath,bm}
\usepackage{longtable}
\usepackage{overpic}

\title{FAST AND ACCURATE NUMERICAL SIMULATIONS FOR THE STUDY OF CORONARY ARTERY BYPASS GRAFTS BY ARTIFICIAL NEURAL NETWORK}

\author{
 Pierfrancesco Siena \\
  MathLab, Mathematics Area, \\
  SISSA International School for Advanced Studies, \\
    Via Bonomea, 265, 34136, Trieste, Italy,\\
  \texttt{psiena@sissa.it} \\
   \And
 Michele Girfoglio \\
  MathLab, Mathematics Area, \\
  SISSA International School for Advanced Studies, \\
    Via Bonomea, 265, 34136, Trieste, Italy,\\
  \texttt{mgifogl@sissa.it} \\
  \And
 Gianluigi Rozza \\
  MathLab, Mathematics Area, \\
  SISSA International School for Advanced Studies, \\
    Via Bonomea, 265, 34136, Trieste, Italy,\\
  \texttt{grozza@sissa.it} \\
}

\begin{document}
\newcommand{\openfoam}{Open\nolinebreak\hspace{-.2em}\nolinebreak\hspace{-.2em}FOAM\textsuperscript{\textregistered}\xspace}
\maketitle
\section*{ABBREVIATIONS}
{\renewcommand\arraystretch{0.8}
\begin{longtable}{p{2cm}p{11cm}p{10cm}}
\hspace{-45pt}	\textbf{ANN}  & Artificial neural network\\
\hspace{-45pt}	\textbf{CABG}  & Coronary artery bypass graft \\
\hspace{-45pt}	\textbf{DEIM}  & Discrete empirical interpolation method\\
\hspace{-45pt}	\textbf{FE}  & Finite element\\
\hspace{-45pt}	\textbf{FFD}  & Free form deformation\\
\hspace{-45pt}	\textbf{FOM}  & Full order model\\
\hspace{-45pt}	\textbf{FV}  & Finite volume \\
\hspace{-45pt}	\textbf{LAD}  & Left anterior descending artery\\
\hspace{-45pt}	\textbf{LCx}  & Left circumflex artery\\
\hspace{-45pt}	\textbf{LITA}  & Left internal thoracic artery\\
\hspace{-45pt}	\textbf{LMCA}  & Left main coronary artery \\
\hspace{-45pt}	\textbf{NURBS}  & Non-uniform rational basis spline\\
\hspace{-45pt}	\textbf{POD}  & Proper orthogonal decomposition\\
\hspace{-45pt}	\textbf{RBF}  & Radial basis functions\\
\hspace{-45pt}	\textbf{ROM}  & Reduced order model\\
\hspace{-45pt}	\textbf{SV}  & Saphenous vein\\
\hspace{-45pt}	\textbf{WSS}  & Wall shear stress\\

\end{longtable}}

\section{INTRODUCTION}\label{sec1}

\captionsetup[subfigure]{position=top,singlelinecheck=off,justification=raggedright}
Coronary artery diseases represent one of the most common causes of mortality worldwide. When coronary arteries are completely or partially occluded for the presence of stenosis, the lack of blood to the cardiac tissue can lead to heart attack. Nowadays, the most successful surgical procedure consists of creating an alternative path to bypass the stenosis, the so-called CABG. In literature, different CABG configurations in presence of single or multiple stenosis are analysed in order to establish a good clinical treatment: see, e.g., \cite{d2002isolated, harling2012surgical, rosenblum2019priorities, scott2000isolated, gaudino2014chronic}. 
Since after some years from the surgical treatment blood supply fails again by causing the need of reintervention, computational investigation of haemodynamic patterns near stenosis and anastomosis regions is of remarkable clinical interest. 

Our study focuses on the numerical simulation of the blood flow features in a patient-specific coronary system when an isolated stenosis of LMCA occurs. A CABG performed with the LITA on the LAD is analyzed. 
High fidelity simulations are able to provide accurate predictions of the blood flow. However in the clinical context repeated model evaluations at varying of physical and geometrical parameters are often required and this leads to very high computational cost when using standard discretization techniques, here referred to as FOMs, e.g. FE or FV techniques. ROMs have been proposed as an efficient tool to approximate FOM systems by significantly reducing
the computational cost required to obtain numerical solutions in a parametric setting (see, e.g., \citealp{benner2015survey}, \citealp{hesthaven2016certified}). A non-intrusive data-driven ROM is used in this work in order to enable fast and reliable computations, involving POD for the computation of reduced basis space and ANNs for the evaluation of the modal coefficients (\citealp{chen2021physics,hesthaven2018non,shah2021finite,pichi2021artificial}). Both primal variables (pressure and velocity) and derived quantities (wall shear stress) are considered.
The introduction of machine learning techniques is appealing due to its development and diffusion in several areas including hemodynamic applications: see, e.g.,  \citealp{liang2020feasibility,gharleghi2020deep,su2020generating,kissas2020machine}. 
The FV method is employed due to its reliable application in some recent works (\citealp{girfoglio2020non,girfoglio2021non,pandey2020review,marsden2015multiscale,buoso2019reduced}) and its widespreading in commercial codes. 

The chapter is organized as follows. The FOM is introduced and discussed in { Section} \ref{sec:background}. Some relevant ROM studies related to similar problems to ours as well as 
the ROM approach here adopted are discussed in { Section} \ref{rom}. The numerical results are shown in { Section} \ref{results}. Finally, conclusions and perspectives are provided in { Section} \ref{concluison}.

\section{THE FULL ORDER MODEL}
\label{sec:background}

\subsection{The Navier-Stokes equations}
\label{Navier-Stokes}

Let us refer to a patient-specific domain $\Omega \subset \mathbb{R}^3$ over a cardiac cycle $(0,T]$ when the transient effects are passed. In this work, the blood is considered as a Newtonian fluid. 
Then if $\mathbf{u} : \Omega \times (0,T] \mapsto \mathbb{R}^3$ and $P : \Omega \times (0,T] \mapsto \mathbb{R}$ are the velocity and the kinematic pressure of the fluid (i.e. the pressure divided by the density), and $\nu \in \mathbb{R}^{+}$ is the kinematic viscosity, the dynamics of the blood is described by the incompressible Navier-Stokes equations:
\begin{equation}
\begin{cases} 
\frac{\partial \bm u}{\partial t}+\nabla \cdot (\bm u \otimes \bm u)-\nu \Delta \bm u + \nabla P= 0  & \quad \text{in} \quad \Omega\times (0,T], \\
\nabla \cdot \bm u =0 & \quad \text{in} \quad \Omega\times (0,T].	
\end{cases}
\label{N-S}
\end{equation}
We also introduce the WSS defined as follows:
\begin{equation}
\text{WSS} = \tau \cdot \mathbf{n},
\end{equation}
where $\tau = \nu (\nabla \mathbf{u} + \nabla \mathbf{u}^T)$ is the stress tensor and $\mathbf{n}$ is the unit normal outward vector.

\begin{figure}
	\centering
	\subfloat[][\label{CB:a}]{\includegraphics[width=.45\textwidth]{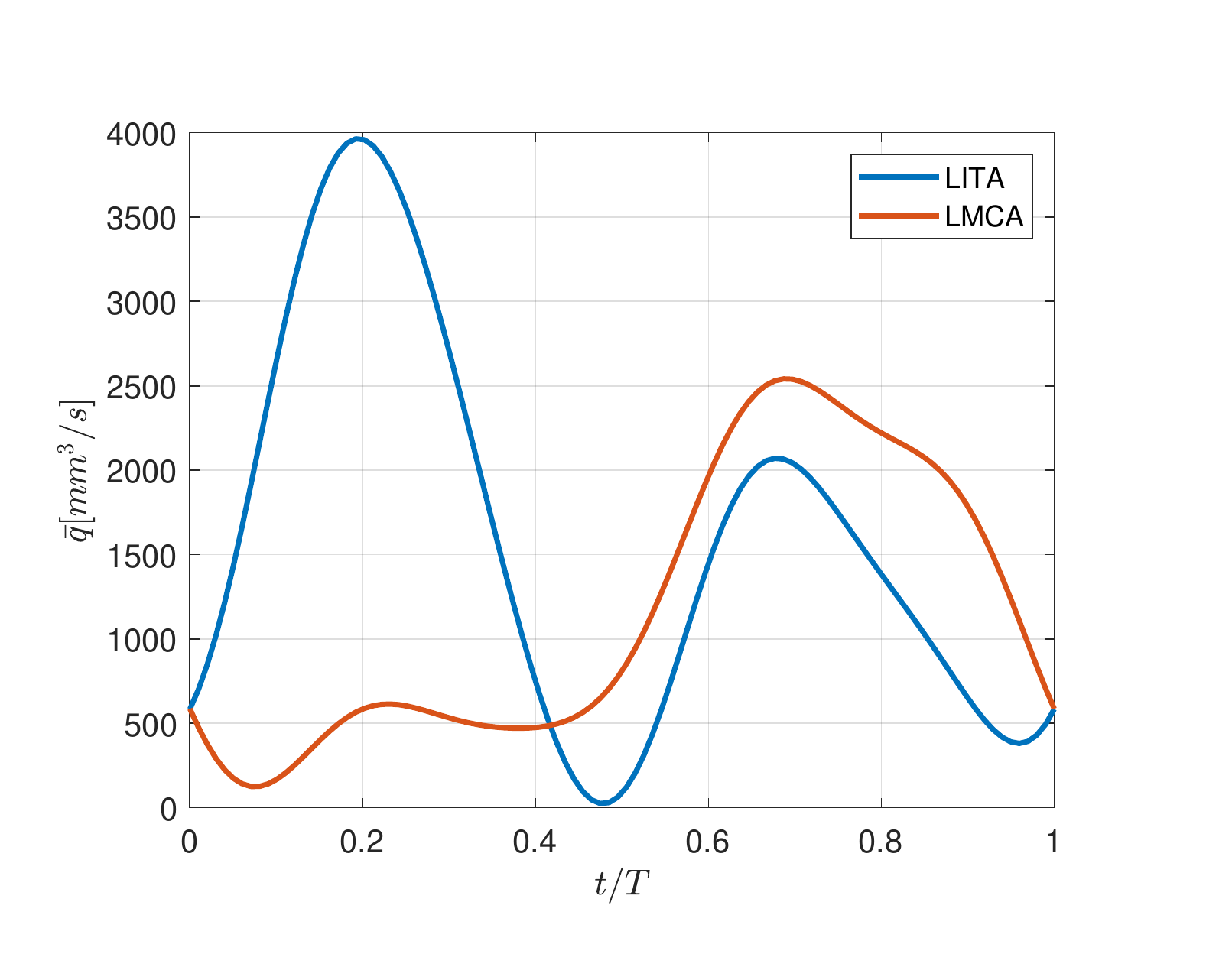}}
	\subfloat[][\label{CB:b}]{\includegraphics[width=.45\textwidth]{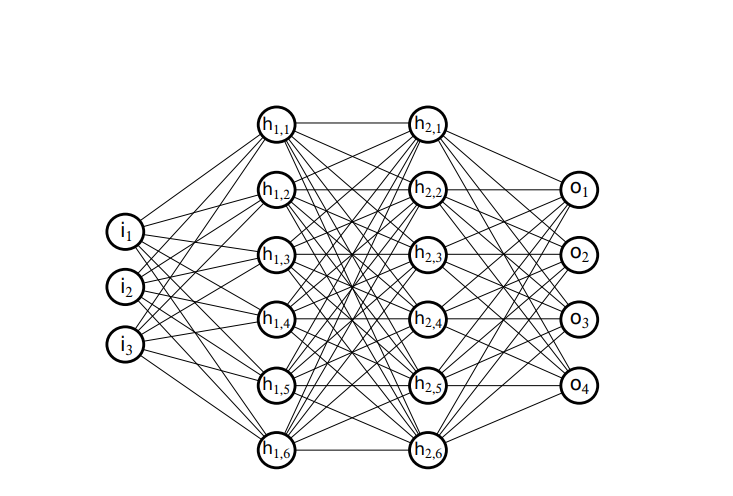}}

	\caption{Time evolution over the cardiac cycle of $\bar{q}_{LMCA}$ and $\bar{q}_{LITA}$ (\citealp{keegan2004spiral,ishida2001mr}) \protect\subref{CB:a} and topology of the feedforward neural network \protect\subref{CB:b}.}
	\label{CB}
\end{figure}

Concerning boundary conditions, a realistic flow rate waveform was enforced on LITA and LMCA sections (see Figure \ref{CB:a}): 
\begin{equation}
q_i(t)=f^i\bar{q}_i(t), \quad i=LMCA,LITA,
\label{flowrate}
\end{equation}
where the factor $f^i \in \left[ \frac{2}{3}, \frac{4}{3} \right] $ is adapted from \cite{keegan2004spiral,ishida2001mr,verim2015cross}.
For the vessel wall, no slip condition is imposed and, for the outflow boundaries, LAD and LCx, homogeneous Neumann boundary
conditions are prescribed (\citealp{ballarin2016fast,ballarin2017numerical}).

\subsection{Time and space discretization}
\label{time-space-discr}


Let us consider a time step $\Delta t \in \mathbb{R}^+$ such that $t^n=n \Delta t$, $n=0, \ldots, N_T$ and $T = N_T\Delta t$. Moreover, let $(\bm{u}^n,P^n)$ be the approximations of velocity and pressure at time $t^n$. 

For the time discretization of {Eqs.} \eqref{N-S}, we adopt a Backward Differentiation Formula of order 2: 

\begin{equation}
\begin{split}
\begin{cases} 
\frac{3}{2 \Delta t} \bm{u}^{n+1}+\nabla \cdot (\bm u^{n} \otimes \bm u^{n+1})-\nu \Delta \bm u ^{n+1}+ \nabla P^{n+1}= \bm{b}^{n+1},  &  \\
\nabla \cdot \bm u^{n+1} =0. & 	
\end{cases} 
\label{N-S-time-discretized}
\end{split}
\end{equation}
where $\bm{b}^{n+1}=\frac{4\bm{u}^n-\bm{u}^{n-1}}{2\Delta t}$. 

For the space discretization of problem \eqref{N-S-time-discretized} we adopt a FV method. The Gauss-divergence theorem allows to write the integral form of the momentum equation in each control volume $\Omega_i$ as follows:
\begin{equation}
\begin{split}
\frac{3}{2\Delta t}\int_{\Omega_i} \bm{u}^{n+1} d \Omega  + & \int_{\partial \Omega_i} (\bm u^{n} \otimes \bm u^{n+1}) \cdot d \bm{A} -  \nu \int_{\partial \Omega_i}  \nabla \bm u^{n+1} \cdot d \bm{A} \\ 
& + \int_{\partial \Omega_i} P^{n+1} d \bm{A}  = \int_{ \Omega_i} \bm{b}^{n+1} d \Omega.
\end{split}
\label{integral_form}
\end{equation}
The convective and diffusive terms are treated using a second order central scheme. If $\bm{u}_i^{n+1}$ and $\bm{b}_i^{n+1}$ indicate the average velocity and the source term in the control volume $\Omega_i$,  $\bm{u}_{i,j}^{n+1}$ and $P_{i,j}^{n+1}$ the velocity and pressure associated to the centroid of face $j$ normalized by the volume of $\Omega_i$, the discretized form of the {Eq.} \eqref{integral_form} can be written as:

\begin{equation}
\frac{3}{2 \Delta t} \bm{u}_i^{n+1}+\sum_j \phi_j \bm{u}_{i,j}^{n+1}-\nu\sum_j(\nabla \bm{u}_i^{n+1})_j\cdot \bm{A}_j+\sum_j P_{i,j}^{n+1}\bm{A}_j=\bm{b}_i^{n+1},
\label{FV_form}
\end{equation}
where $\bm A_j$ is the surface vector of each face $j$ of the control volume and $\phi_j = \bm{u}_j^{n} \cdot \bm{A}_j$ is the convective flux associated to $\bm{u}^{n}$ through face $j$ of the control volume. Applying the divergence operator to {Eq.} \eqref{FV_form} in semi-discretized form, i.e. with the pressure term in continuous form while all the other terms in discrete form, and introducing the continuity equation, the following Poisson equation for the pressure is derived:
\begin{equation}
\Delta P^{n+1} = \nabla \cdot \Big( -\sum_j \phi_j \bm{u}_{j}^{n+1} + \nu \sum_j(\nabla \bm{u}^{n+1})_j\cdot \bm{A}_j + \bm{b}^{n+1} \Big).
\label{pressure_eq}
\end{equation}
Integrating {Eq.} \eqref{pressure_eq} and applying the Gauss-divergence theorem, we obtain:
\begin{equation}
\sum_j (\nabla P^{n+1})_j\cdot \bm{A}_j = \sum_j  \Big( -\sum_j \phi_j \bm{u}_{j}^{n+1} + \nu \sum_j(\nabla \bm{u}^{n+1})_j\cdot \bm{A}_j + \bm{b}^{n+1} \Big)_j \cdot \bm{A_j}.
\label{FV_form1}
\end{equation}
We use the Pressure-Implicit with Splitting of Operators 
algorithm (\citealp{issa1986solution}) available in \openfoam
to deal with the coupled set of {Eqs.} \eqref{FV_form} and \eqref{FV_form1}. 

\subsection{The computational domain}
\label{domain}
\begin{figure}
	\centering
	\subfloat[][\label{CABG1:a}]{\includegraphics[width=.3\textwidth]{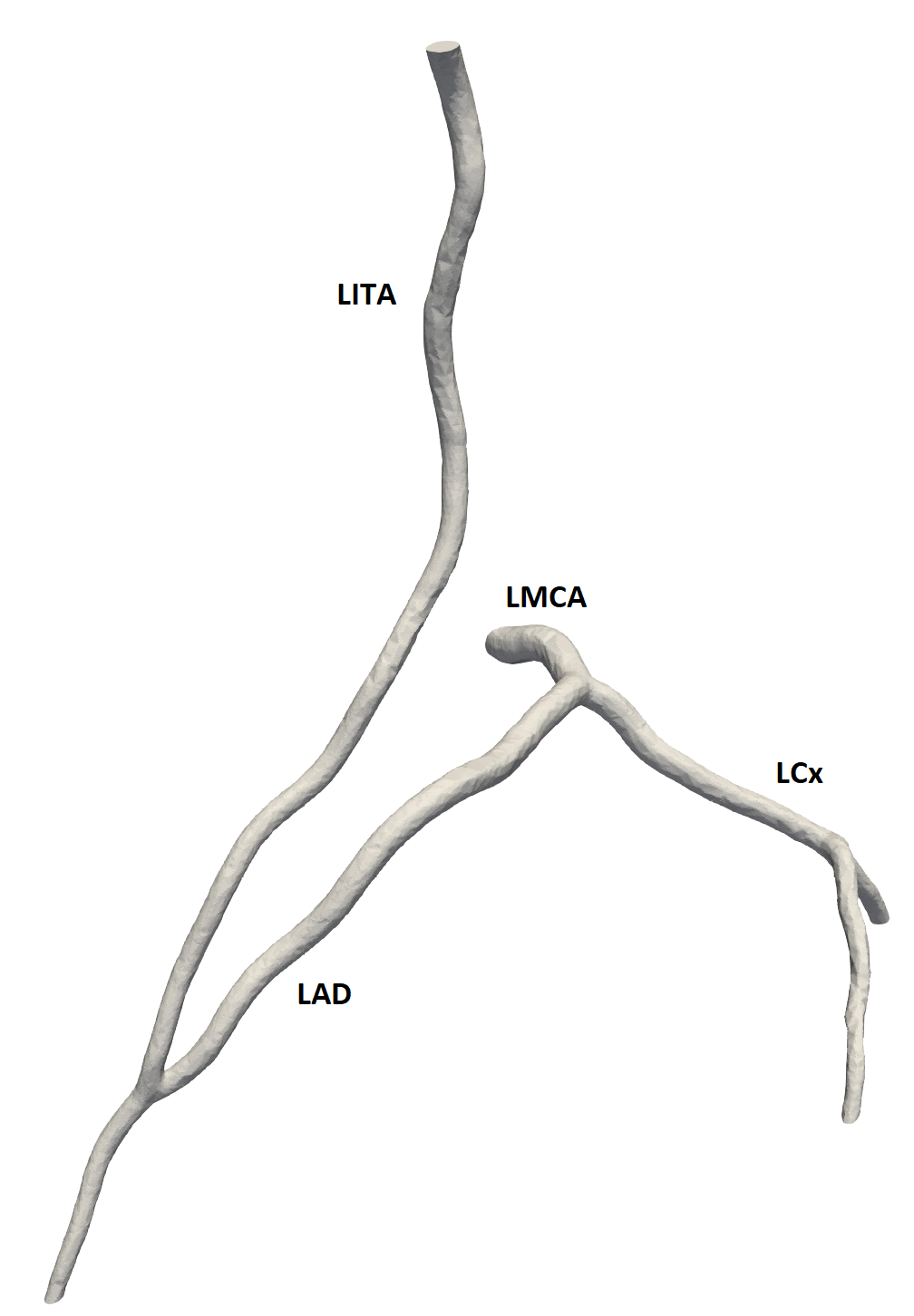}} \qquad
	\subfloat[][\label{CABG1:b}]{\includegraphics[width=.4\textwidth]{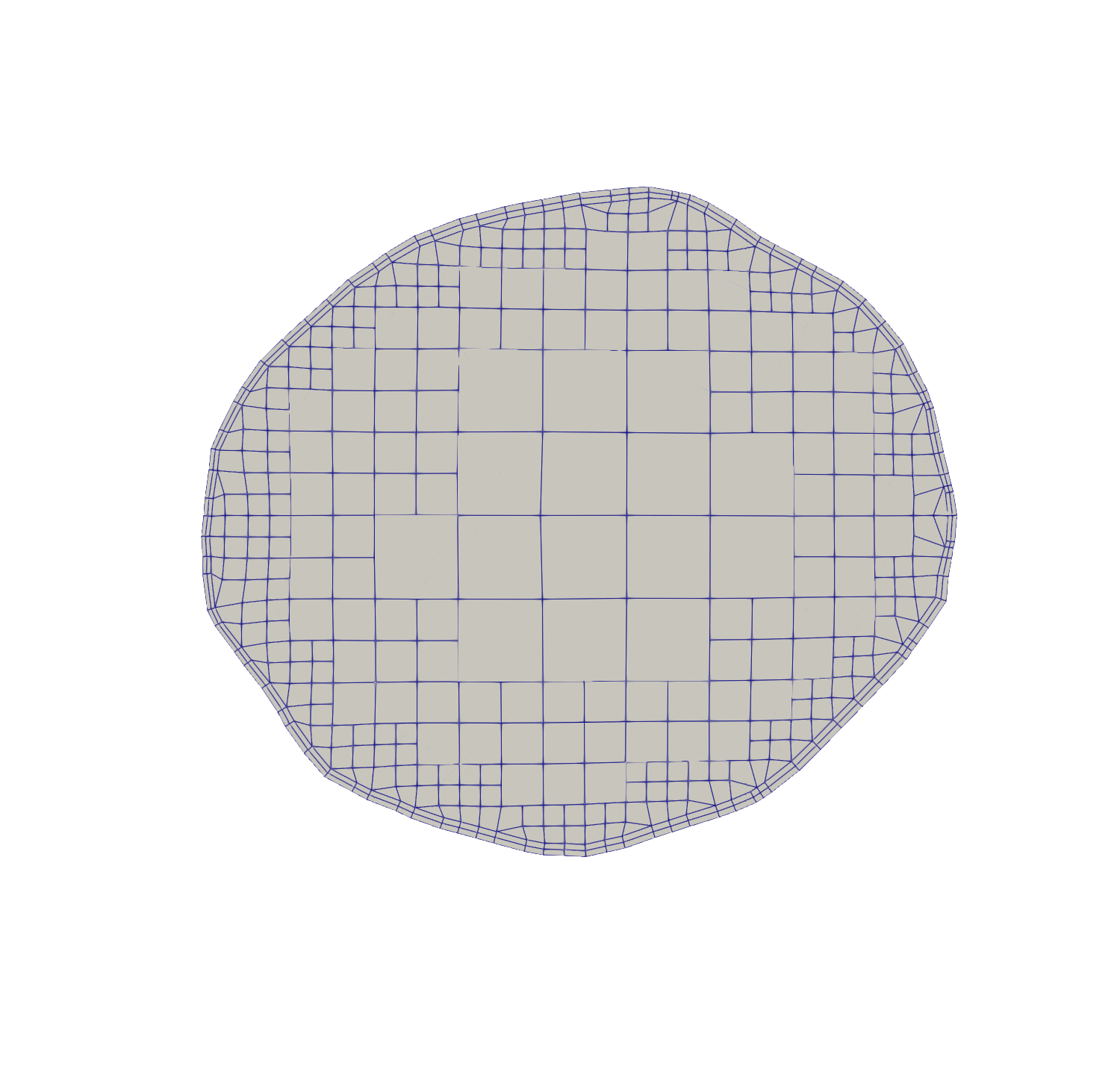}}
	\caption{CABG geometry \protect\subref{CABG1:a} and view of the mesh on an internal section \protect\subref{CABG1:b}.}
	\label{CABG1}       
\end{figure}

Computed tomography images of a patient specific configuration 
were provided by Ospedale Luigi Sacco in Milan. The process to construct the virtual geometry (see {Fig.} \ref{CABG1:a}) was described in detail in \cite{ballarin2016fast,ballarin2017numerical}.
The mesh ({Fig.} \ref{CABG1:b}) is built using the mesh generation utility \emph{snappyHexMesh} available in \openfoam. It has 986.278 cells and its minimum and maximum diameter are $h_{min} \simeq 3.0 \cdot 10^{-5} m$ and $h_{max} \simeq 4.3 \cdot 10^{-4} m$, respectively. The quality of the mesh is rather high: it has very low values of 
average non-orthogonality (13$^\circ$), skewness ($\leq 3$), and aspect ratio ($\leq 16$).

In order to introduce the stenosis in the LMCA, FFD is performed by means of a NURBS volumetric parameterization (\citealp{brujic2002measurement,tezzele2021pygem}). The shape parametrization of the computational domain 
is deeply analysed in literature (\citealp{manzoni2012model,lassila2011reduction,manzoni2014reduced,lassila2013boundary,lassila2013reduced,manzoni2012shape,stabile2020efficient,burgos2015nurbs,tezzele2018combined}) using both FFD and RBF techniques. 
However, whilst the RBF approach deforms the grid as a whole and it does not preserve the original geometry, in the FFD method all the vertices remain on the initial surface. 

\begin{figure}
	\centering
	\subfloat[][\label{mimmo:a}]{\includegraphics[width=.4\textwidth]{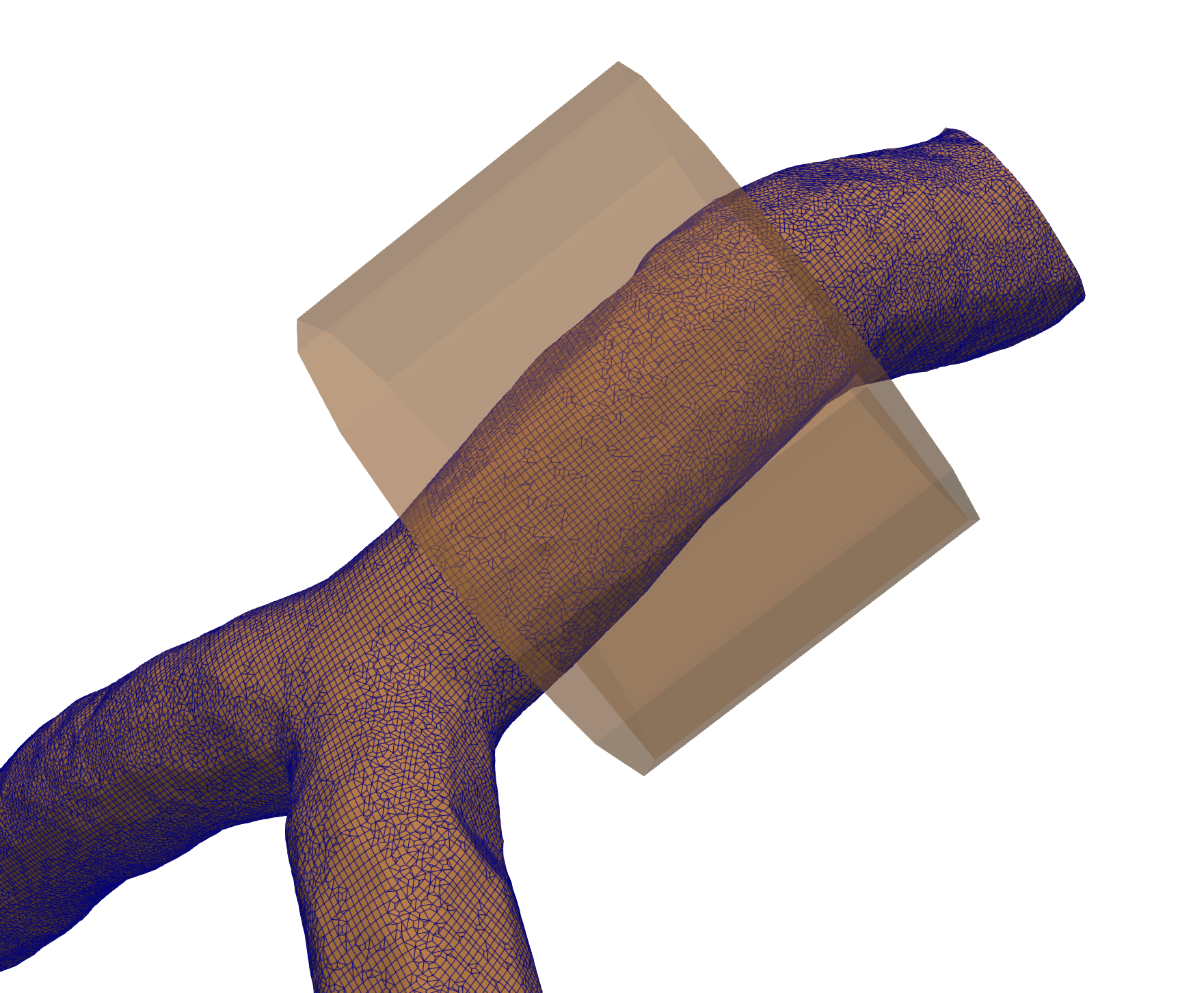}}
	\subfloat[][\label{mimmo:b}]{\includegraphics[width=.4\textwidth]{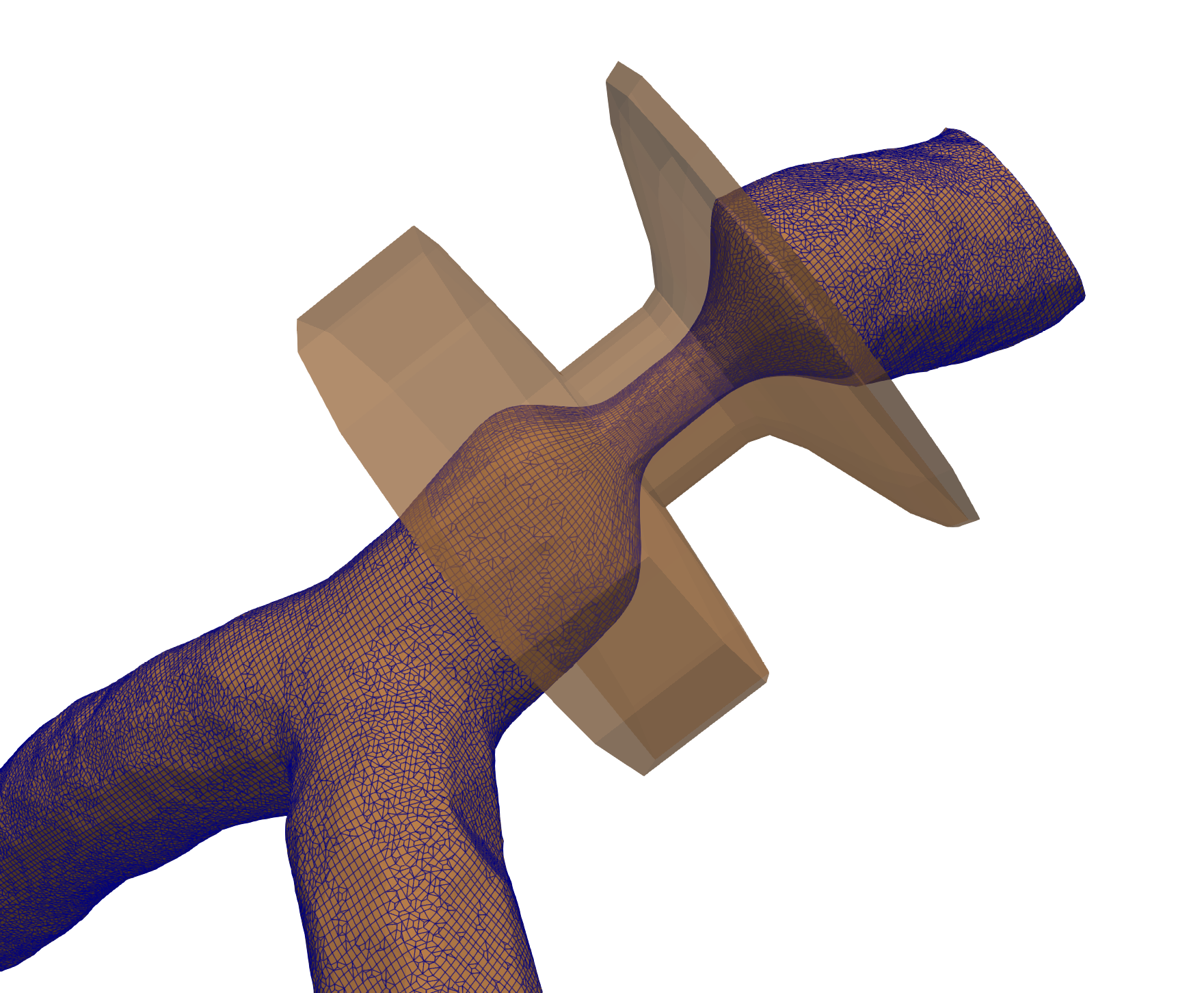}}
	\caption{Introduction of $70\%$ stenosis in the LMCA using FFD technique: \protect\subref{mimmo:a} initial lattice, \protect\subref{mimmo:b} deformed lattice.}
	\label{Mimmo}
\end{figure}

The FFD method consists of three steps (\citealp{amoiralis2008freeform,lamousin1994nurbs}): 
\begin{itemize}
	\item[(i)] A parametric lattice of control points 
	is constructed generating a structured mesh 
	around the region of the LMCA where the stenosis occurs.
	Then the control points are used to define a NURBS volume which contains the domain to be warped. Figure \ref{mimmo:a} displays the lattice in its initial configuration. 

	\item[(ii)] 
	The octree algorithm (\citealp{amoiralis2008freeform}) could be used to find a match between the control points of the lattice and the points of the computational domain. 

	\item[(iii)] The coordinates of the control points are modified, so that the parametric volume and consequently the computational domain are deformed. 
	Figure \ref{mimmo:b} shows the lattice in its deformed configuration where a $70\%$ stenosis is introduced. 
\end{itemize}

We highlight that the introduction of a $70\%$ stenosis in the LMCA does not affect the mesh quality. 

\section{REDUCED ORDER MODEL}
\label{rom}

\subsection{POD-ANN approach}

In this work we use the POD-ANN approach consisting of two main stages:
\begin{itemize}
	\item \emph{Offline}: a reduced basis space is built by applying POD to a database of high-fidelity solutions obtained by solving the FOM for different values of physical and/or geometrical
parameters. Once the reduced basis space is computed, we project the original snapshots onto such a space by obtaining the corresponding parameter dependent modal coefficients. Then the training of the neural networks to approximate the map between parameters and modal coefficients is carried out. This stage is computationally expensive, however it only needs to be performed once. 
	\item \emph{Online}: for any new parameter value, we approximate the new coefficients by using the trained neural network and the reduced solution is obtained as a linear combination of the POD basis functions multiplied by modal coefficients. During this stage, it is possible to explore the parameter space at a significantly reduced cost.
\end{itemize}
Here we are going to test our ROM for the reconstruction of the time evolution of the system. So our parameter of interest coincides with the time. Other parametric studies, involving geometrical features such as the stenosis degree or physical properties such as the scaling factor $f^i$ (see {Eq.} \eqref{flowrate}), will be addressed in a work in preparation (\citealp{pier2022neural}). 

\subsubsection{Proper orthogonal decomposition and method of snapshots}

We solve the FOM described in { Section} \ref{time-space-discr} for each time $t_k \in \{t_1,\dots, t_N \} \subset (0,T]$. 
The high fidelity solution $\Phi_h = \{\mathbf{u}, P, \text{WSS}\}$ can be stored into the matrix $\mathcal{S}_\Phi$ in a column-wise sense:
\begin{equation}
	\mathcal{S}_\Phi=\{  \Phi_h(t_1)| \dots | \Phi_h(t_N) \} \in\mathbb{R}^{N_{\phi_h} \times N},
\end{equation}
where the subscript $h$ denotes a solution computed with the FOM and $N_{\phi_h}$ is the dimension of the space field $\Phi$ belong to in the FOM.
If $R \leq min(N_{\phi_h},N)$ is the rank of $\mathcal{S}_{\Phi}$, the singular value decomposition enables to factorise $\mathcal{S}_\Phi$ as:
\begin{equation}
	\mathcal{S}_\Phi=\mathcal{W} \mathcal{D} \mathcal{Z}^T,
\end{equation}
where $\mathcal{W}=\{ \bm{w}_1| \dots |\bm{w}_{N_{\phi_h}} \} \in \mathbb{R}^{N_{\phi_h} \times N_{\phi_h}}$ and $\mathcal{Z}=\{ \bm{z}_1| \dots |\bm{z}_{N}  \} \in \mathbb{R}^{ N \times N}$ are two orthogonal matrices composed of left singular vectors and right singular vectors respectively, and $\mathcal{D} \in \mathbb{R}^{ N_{\phi_h} \times N}$ is a diagonal matrix with $R$ non-zero real singular values $\sigma_1 \geq \sigma_2 \geq \dots \geq \sigma_R > 0$.\\The purpose is to find $L<R$ suitable orthonormal vectors which approximate the columns of $\mathcal{S}_\Phi$. The Schmidt-Eckart-Young theorem states that the POD basis of rank $L$  consists of the first $L$ left singular vectors of $\mathcal{S}_\Phi$, also named modes (\citealp{eckart1936approximation}). In this work, the extrapolated modes are resumed as columns in the matrix $\mathcal{V}$:
\begin{equation}
\mathcal{V}=\{ \bm{w}_1| \dots |\bm{w}_L \} \in \mathbb{R}^{ N_{\phi_h} \times L}.
\end{equation}
It is well known the reduced basis is the set of vectors that minimizes the distance between the snapshots and their projection into the space spanned by the reduced basis. In addition, the error introduced by replacing the columns of $\mathcal{S}_\Phi$ with those of $\mathcal{V}$ is the sum of the squares of the neglected singular values; therefore, tuning $L$, it is possible to approximate $\mathcal{S}_\Phi$ with arbitrary accuracy (\citealp{quarteroni2015reduced}). A common choice is to set $L$ equal to the smallest integer such that:
\begin{equation}
\frac{\sum_{i=1}^{L}\sigma_i}{\sum_{i=1}^{R}\sigma_i} \ge \delta, 
\label{energy}
\end{equation}
where $\delta$ is an user-provided treshold and the left hand side of {Eq.} \eqref{energy} represents the percentage of energy retained by the first $L$ modes. 
 
Once the POD basis is available, the reduced solution $\Phi_{rb} (t_k)$ that approximates the full order solution $\Phi_{h} (t_k)$ is:
\begin{equation}
 \Phi_h (t_k) \simeq \Phi_{rb} (t_k) = \sum_{j=1}^{L} (\mathcal{V}^T  \Phi_h (t_k))_j \mathbf w_{j}, \quad t_k \in \{t_1,\dots, t_N \} \subset (0,T],
\end{equation}
where $(\mathcal{V}^T  \Phi_h(t_k))_j$ is the modal coefficient associated to the $j$-th mode.
\subsubsection{Interpolation with neural networks}

An artificial neural network is a computational model composed of neurons and synapses. 
It can be represented with an oriented graph, with the neurons as nodes and the synapses as oriented edges. In this work we use feedforward network (also named perceptron), where input ($i$), hidden ($h$) and output ($o$) nodes are arranged into layers ({Fig.} \ref{CB:b}). 
A training process is carried out during the \emph{offline} stage to adjust the synaptic weights and to configure the network: 
in particular, 
parameters such as the activation function, the number of layers, the number of neurons per layer and the learning rate are tuned to minimize the mean square loss function. Further details can be found in literature (\citealp{sharma2017activation,hesthaven2018non,rumelhart1986learning,montana1989training,goodfellow2016deep,chen2021physics}).

The aim is to train the neural network to find the approximation $\bm{\pi}$ of the function 
\begin{equation}
\bm{\pi}: t_k \in \{t_1,\dots, t_N \} \subset (0,T] \mapsto [\mathcal{V}^T  \Phi_h (t_k)]_{j=1}^{L}. 
\end{equation}
Once the networks are trained, the solution can be estimated for a new time value $t_{new}$ during the \emph{online} stage as: 
\begin{equation}
\Phi_{rb}(t_{new})=\sum_{j=1}^{L} {\pi}_j(t_{new})\mathbf w_{j}
\end{equation}
For the creation and training of the neural networks, we employed the Python library PyTorch.

\subsection{State of the art}
In this section we are going to review some recent ROM works focusing on problems similar to ours by highlighting the main differences with respect to the framework here adopted. For a more general discussion about ROM including remarkable applications in several contexts, we refer to \cite{benner2020model}. 

In \cite{buoso2019reduced} an intrusive ROM based on POD-Galerkin strategy within a FV environment is developed in order to estimate pressure drop along blood vessels. On the other hand, here we adopt a non-intrusive data-driven ROM based on POD-ANN method. Another significant distinction with respect to our approach is the use of an idealized (i.e., non patient-specific) geometry. 
Moreover, while FFD is used in our work, in \cite{buoso2019reduced} DEIM is performed to introduce the stenosis. 

In \cite{zainib2020reduced} patient-specific geometries of CABGs are considerd but, while in our work a FV approach is adopted, the FE method is employed in \cite{zainib2020reduced}. Another difference is related to the mathematical approach: in  \cite{zainib2020reduced} an optimal flow control model is presented to obtain meaningful boundary conditions. In order to optimize the shape of the bypass, the theory of optimal control based on the
analysis of the wall shear stress is used also in \cite{quarteroni2003optimal,rozza2005optimization}. For recent works related to optimal control problems in several contexts, the reader is referred, e.g., to \citealp{fevola2021optimal,demo2021extended,donadini2021data}. 

A POD-Galerkin technique is adopted also in \cite{ballarin2016fast,ballarin2017numerical} coupled with an efficient centerlines-based parametrization for the deformation of patient-specific configurations of CABGs. In addition to the differences associated with the ROM approach and the deformation technique, it should be noted that whilst in our work the FV method is adopted, in \cite{ballarin2016fast,ballarin2017numerical} the FE method is employed. 

In recent times, machine learning seems to support considerably cardiovascular medicine. 
In \cite{su2020generating} machine learning is used as an alternative to computational fluid dynamics for generating hemodynamic parameters in real-time diagnosis during medical examinations. The approach is validated by considering the wall shear stress in coronary arteries. 
Anyway, reduced order models and neural networks can efficiently operate jointly as showed in this work and in \cite{pier2022neural}. Another relevant example is given in \cite{fossan2021machine} where prior physics-based knowledge deriving from a reduced-order model is integrated into a neural network framework at the aim to predict pressure losses across coronary segments. In addition, a recent application for multiphase flows that could be enforced in the cardiovascular framework, is available in \cite{papapicco2021neural}.




\section{NUMERICAL RESULTS}
\label{results}
In this section, we present several numerical results for our ROM approach. 

\begin{figure}
	\centering
	\subfloat[][$2^\circ$ coefficient of pressure\label{coeff:a}]{\includegraphics[width=.33\textwidth]{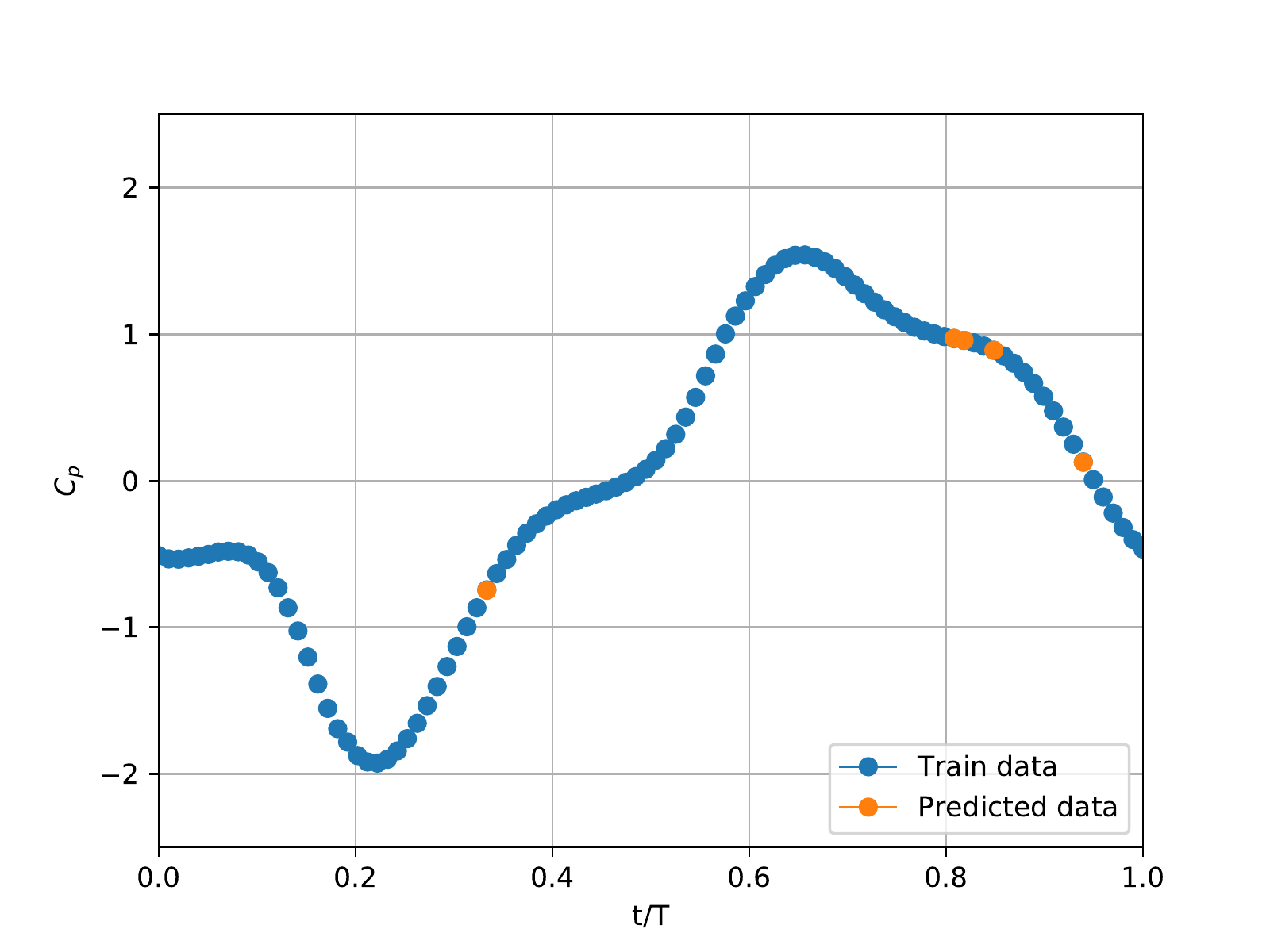}}
	\subfloat[][$1^\circ$ coefficient of velocity\label{coeff:b}]{\includegraphics[width=.35\textwidth]{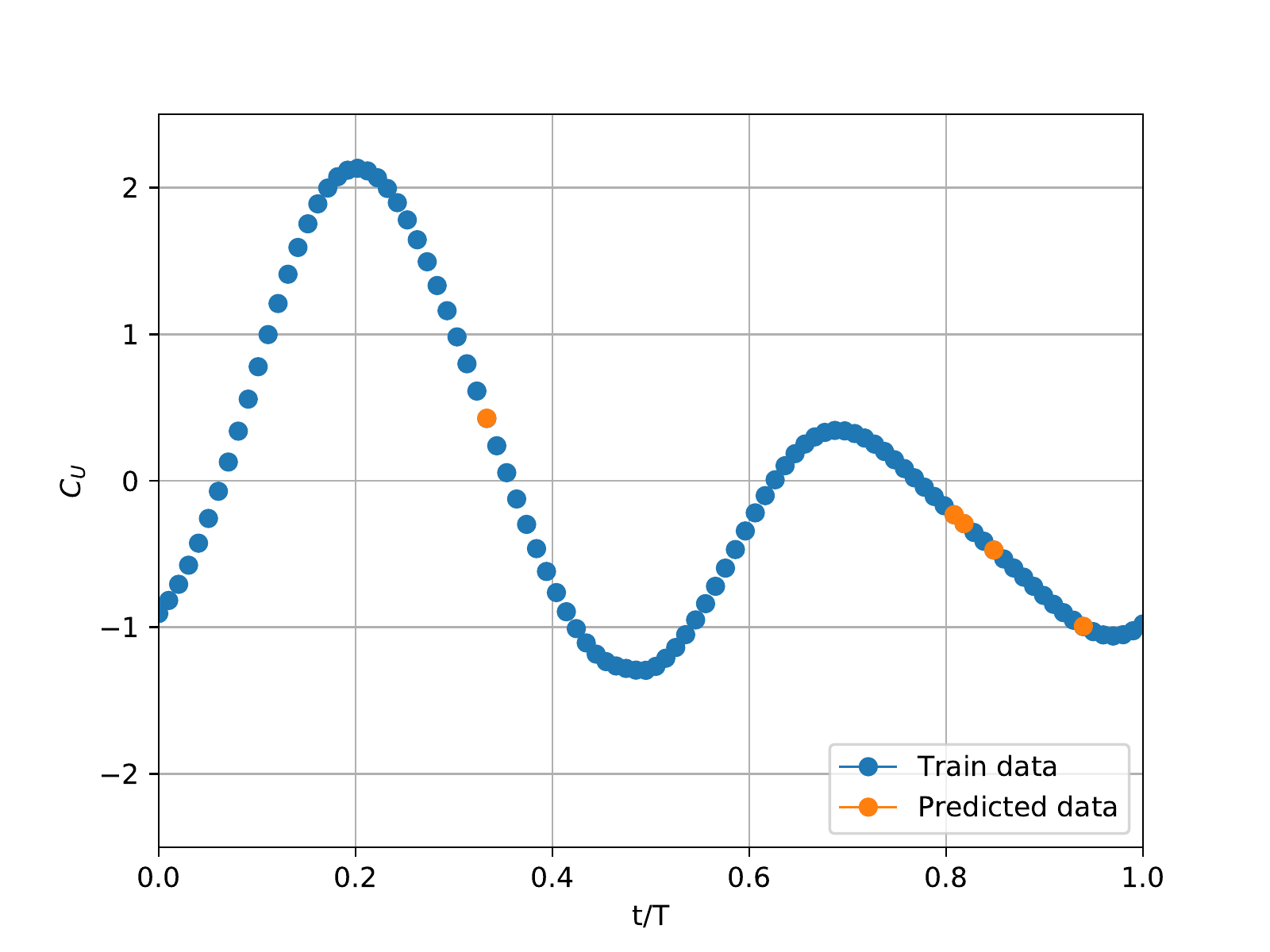}}
	\subfloat[][$4^\circ$ coefficient of WSS\label{coeff:c}]{\includegraphics[width=.35\textwidth]{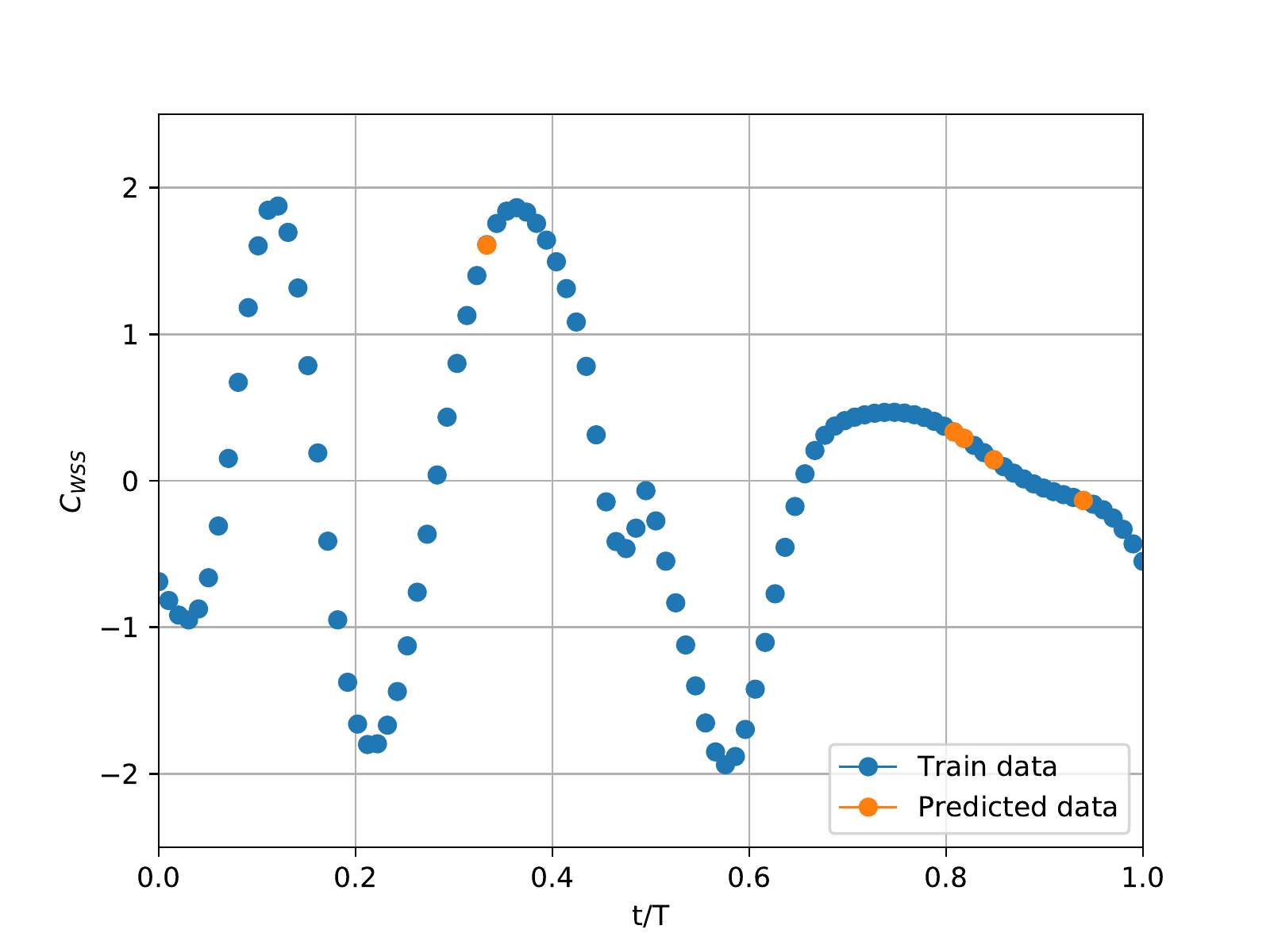}}
	\caption{Time evolution of some reduced coefficients with training points (blue) and test points (orange).}
	\label{coeff}
\end{figure}

We consider the domain in {Fig.} \ref{CABG1:a} where a $70 \% $ stenosis of the LMCA is introduced ({Fig.} \ref{Mimmo}). 
Trial and error process is employed to optimize the hyperparameters of the neural networks. For each variable, { Table} \ref{param_nn} shows hyperparameters which have provided the best performance. The training data are related to the 95\% of the total data provided by the full order model whilst the remaining 5\% is used to do the validation.
\begin{table}
	\caption{Hyperparameters of the neural networks used to reconstruct the reduced coefficients of pressure, velocity and WSS.}
	\begin{tabular}{|c|c|c|c|c|c|}
		\hline
		& \textbf{Neurons per layer} & \textbf{Activation function} & \textbf{Number of epochs} & \textbf{Learning rate} & \textbf{Hidden layers} \\
		\hline
		\textbf{p}  & 500 & ReLU & 50.000 & 1.00e-6 & \\
		\cline{1-5}
		\textbf{U} & 850 & Tanh & 100.000 & 8.25e-6 & 3 \\
		\cline{1-5}
		\textbf{WSS} & 900 & Tanh & 100.000 & 5.50e-6 & \\
		\hline
	\end{tabular}
	\label{param_nn}
\end{table}
With the aim to show the functionality of the networks, some reduced coefficients are displayed in {Fig.} \ref{coeff}, where one can observe that test points (orange) are consistent with training data (blue).

\begin{figure}
	\centering
	\subfloat[][\label{err_snap:a}]{\includegraphics[width=.3\textwidth]{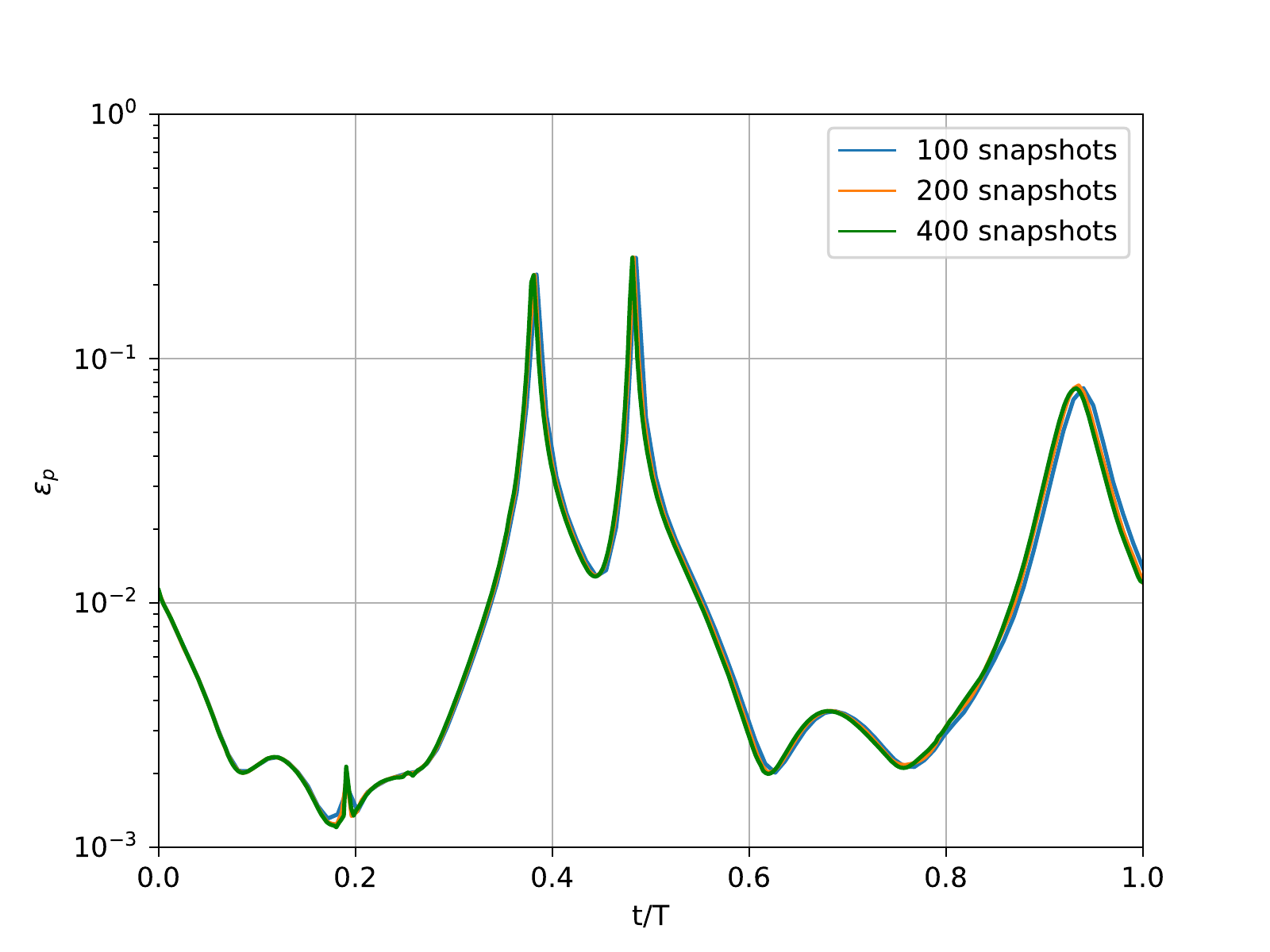}}
	\subfloat[][\label{err_snap:b}]{\includegraphics[width=.32\textwidth]{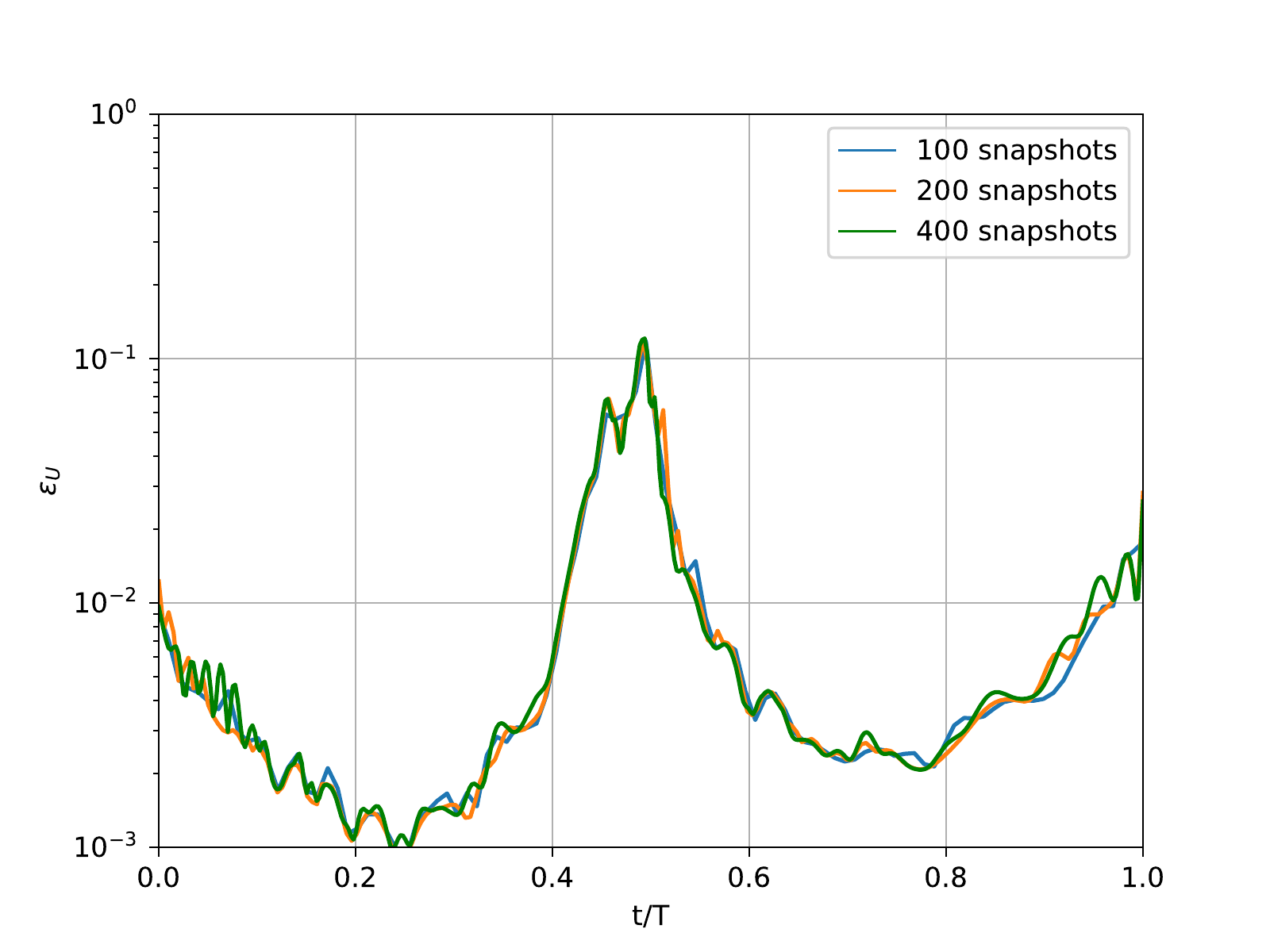}}
	\subfloat[][\label{err_snap:c}]{\includegraphics[width=.32\textwidth]{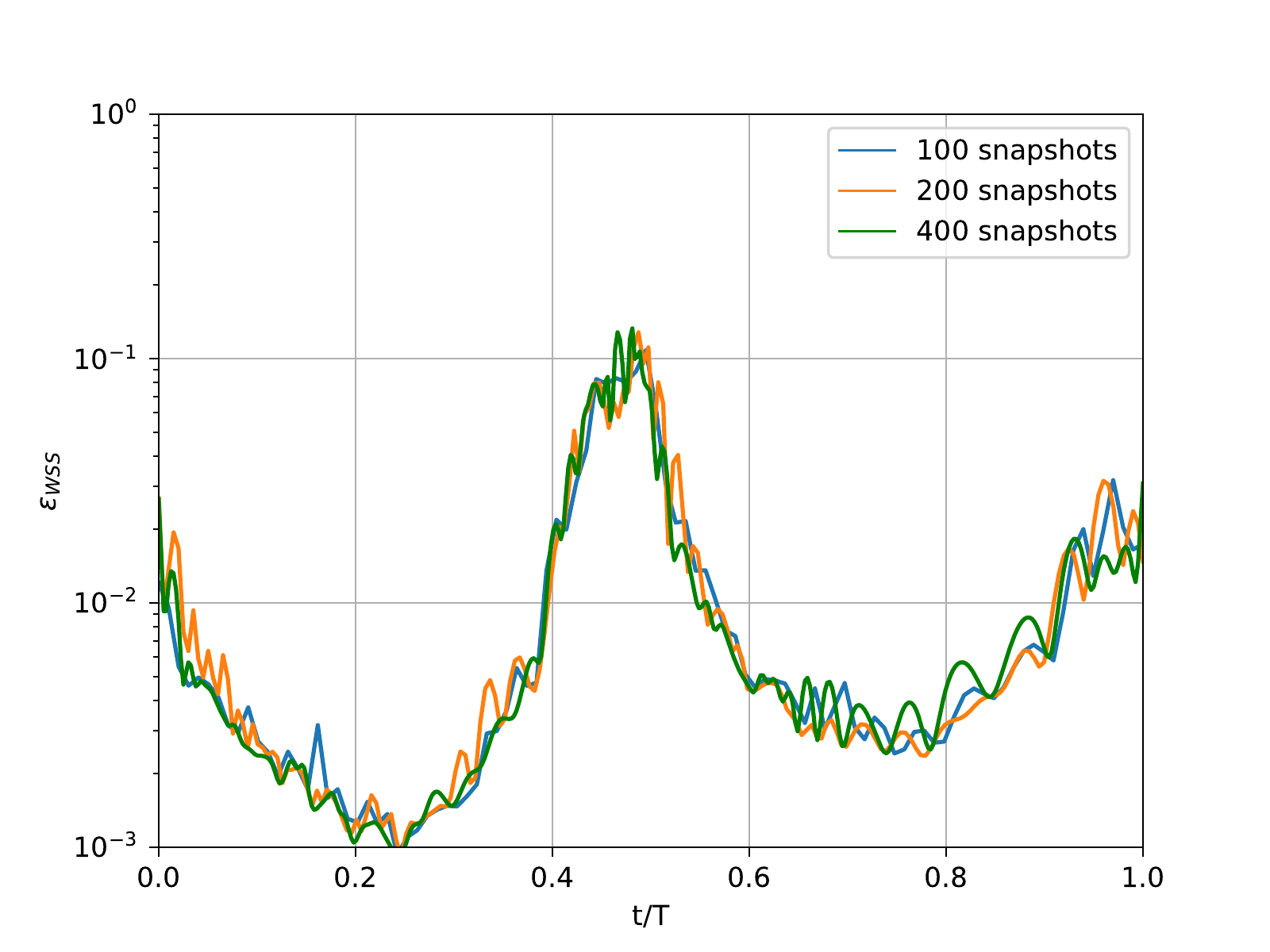}} \\
        \subfloat[][\label{err_mod:d}]{\includegraphics[width=.3\textwidth]{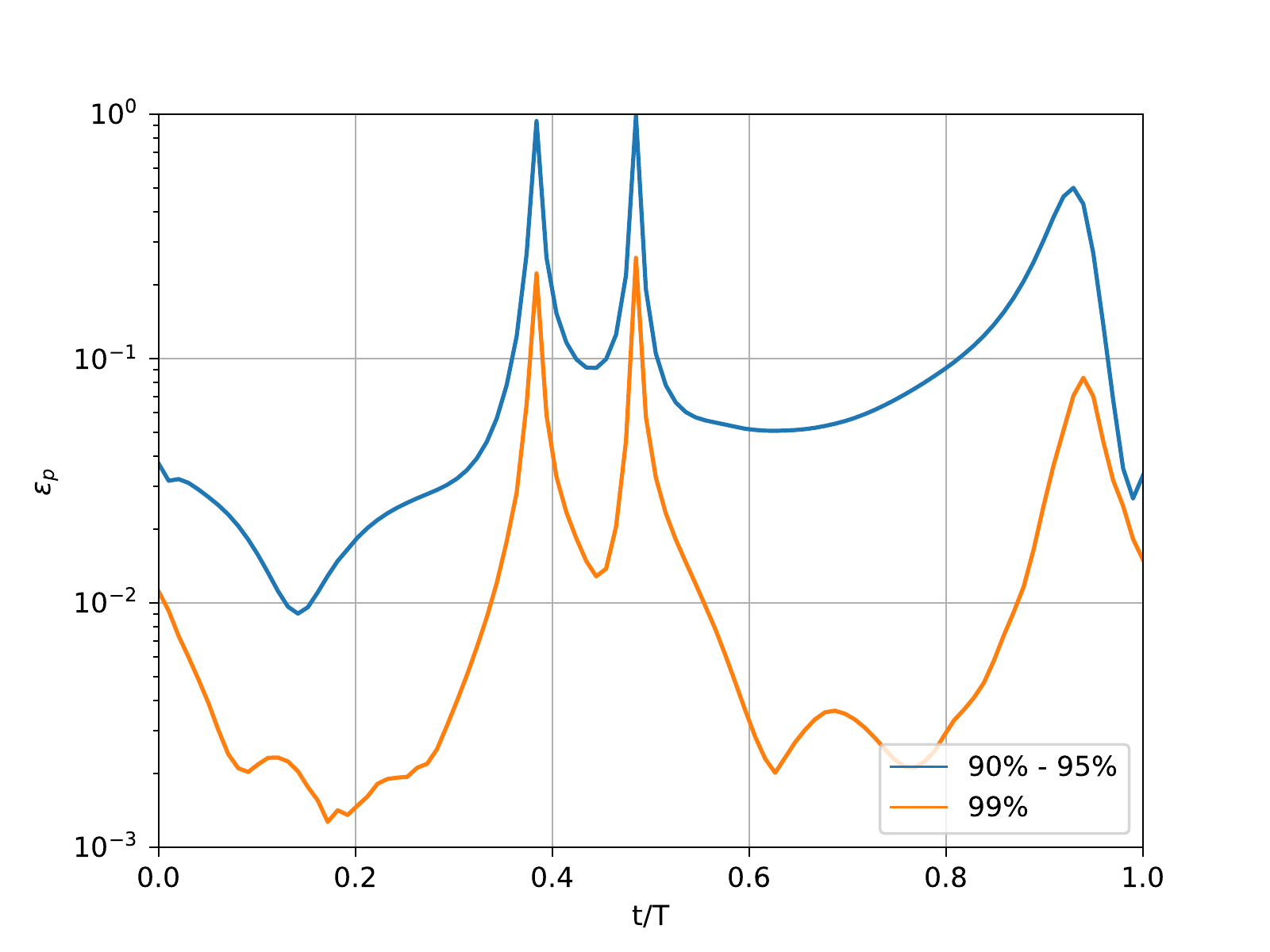}}
	\subfloat[][\label{err_mod:e}]{\includegraphics[width=.32\textwidth]{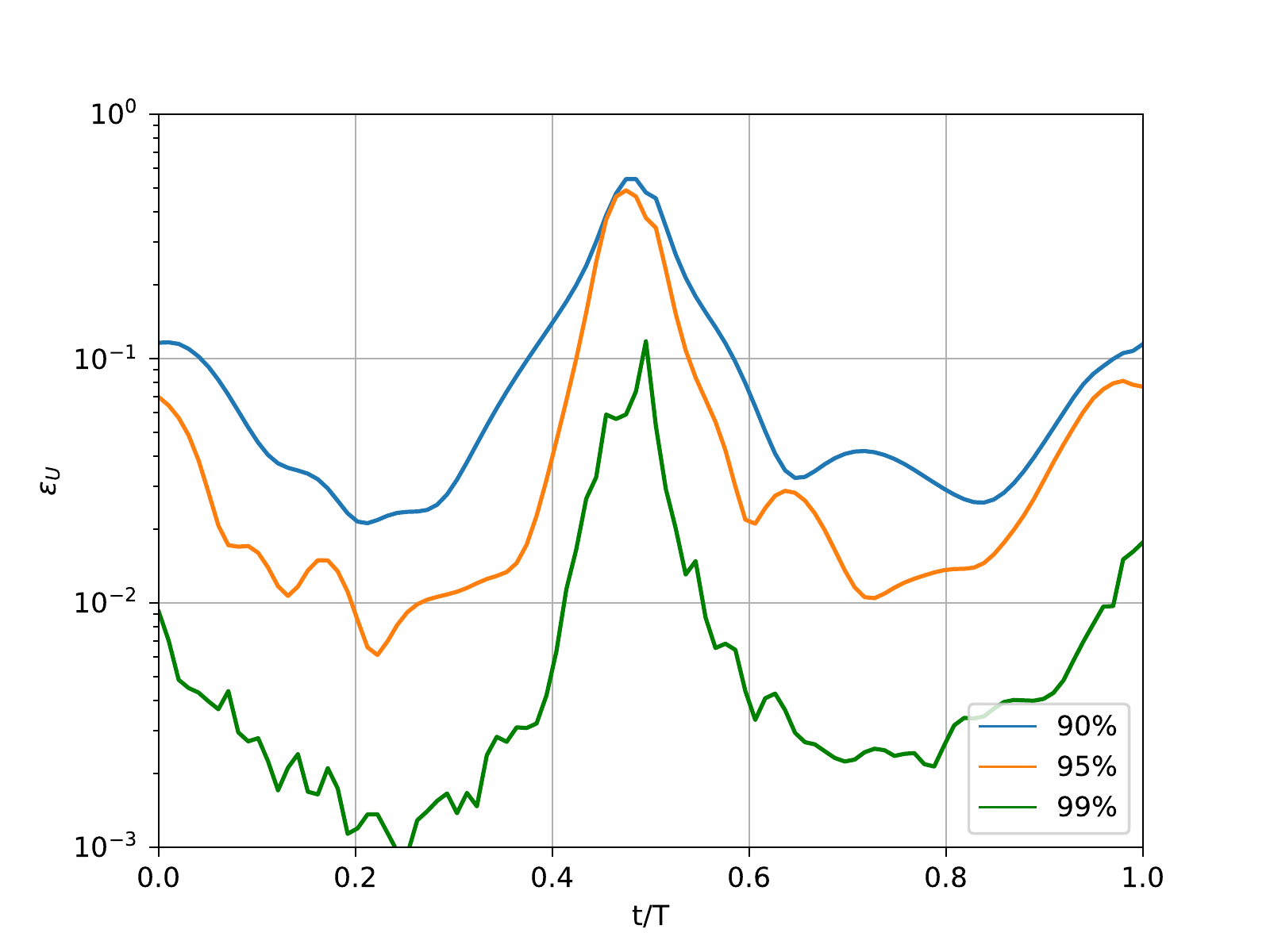}}
	\subfloat[][\label{err_mod:f}]{\includegraphics[width=.32\textwidth]{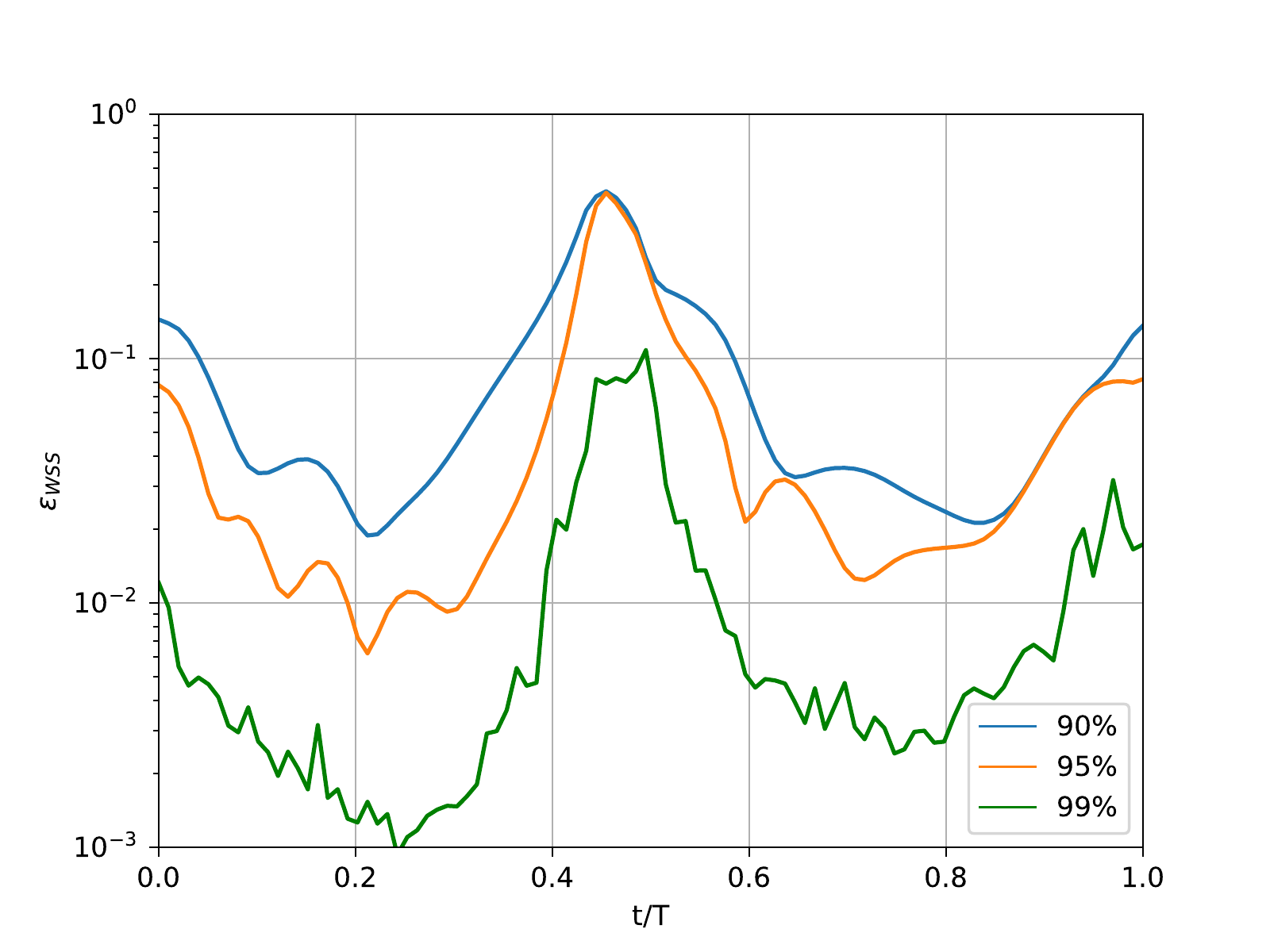}}
	\caption{Time evolution of the relative error $\varepsilon$ varying the number of snapshots (first raw) for \protect\subref{err_snap:a} pressure, \protect\subref{err_snap:b} velocity and \protect\subref{err_snap:c} WSS, and varying the number of modes (second raw) for \protect\subref{err_mod:d} pressure, \protect\subref{err_mod:e} velocity and \protect\subref{err_mod:f} WSS.}
	\label{err_snap}
\end{figure}

We start to perform a convergence test as the number of snapshots increases. 
We collect 100, 200 and 400 full order snapshots, each every 0.008 s, 0.004 s and 0.002 s, respectively (i.e.~we use an equispaced grid in time), over a cardiac cycle $T=0.8$ s. 
{Fig.} \ref{err_snap} shows $L^2$-norm relative error $\varepsilon$ for pressure, velocity and WSS fields over time
for the three different sampling frequencies. We set $\delta = 0.99$ based on the first 90 most energetic POD
modes (3 modes for pressure, 15 modes for velocity and 16 modes for WSS). {Fig.} \ref{err_snap} shows that there is no substantial difference in the errors. Thus, to reduce the computational cost of the offline phase, we will consider 100 snapshots for the results presented from here on. 

Next, we carry out a convergence test as the number of modes increases. We consider three different energy thresholds: $\delta=0.90$ (1 mode for pressure, 3 modes for velocity and 3 modes for WSS), $0.95$ (1 mode for pressure, 5 modes for velocity and 5 modes for WSS), and $0.99$ (3 modes for pressure, 15 modes for velocity and 16 modes for WSS). 
Of course the pressure error is the same for $\delta=0.90,0.95$ because a single mode is enough to reach $95\%$ of the system energy. From {Fig.} \ref{err_snap}, 
we observe that all the variables show a monotonic convergence as the number of the modes is increased. For $\delta=0.99$ we obtain time-averaged errors of about 1.7\% for the pressure and less than 1\% for velocity and WSS. 

\begin{figure}
	\centering	
	
	\subfloat[][]{%
		\begin{overpic}[width=0.32\textwidth]{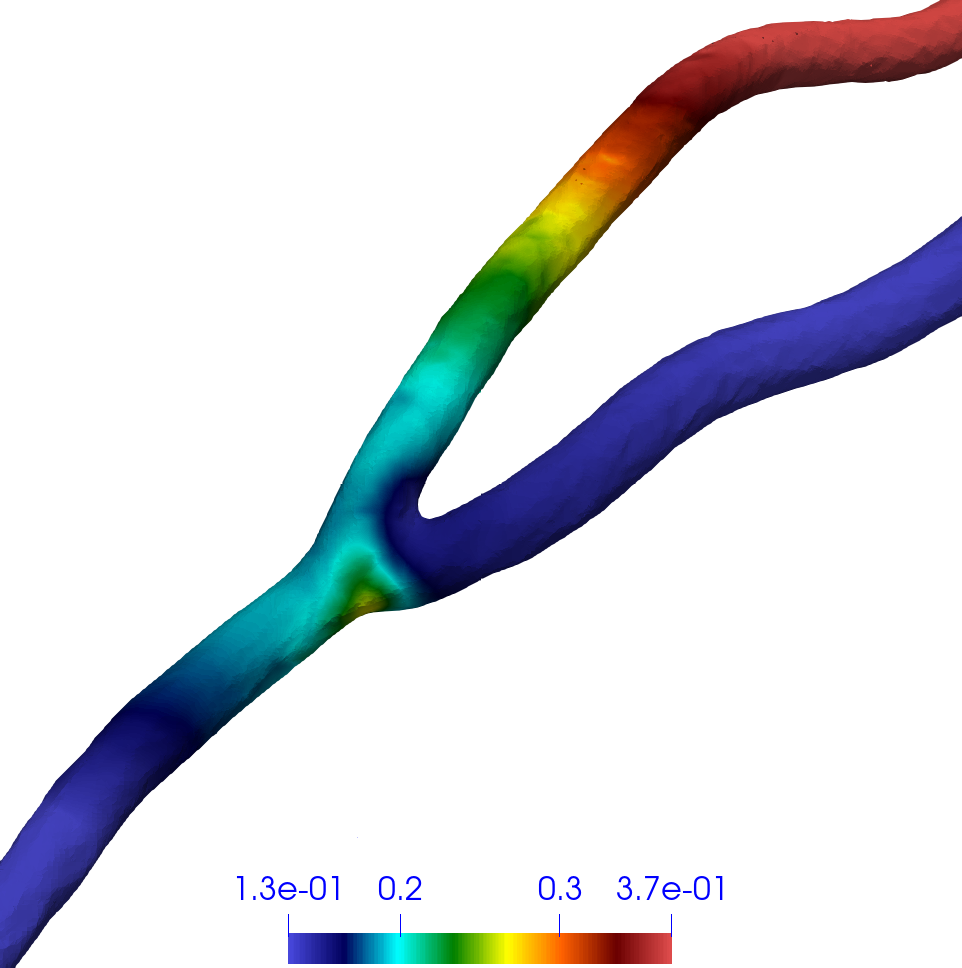}
			\put(20,105){FOM}
		\end{overpic}}
	\subfloat[][]{%
		\begin{overpic}[width=0.32\textwidth]{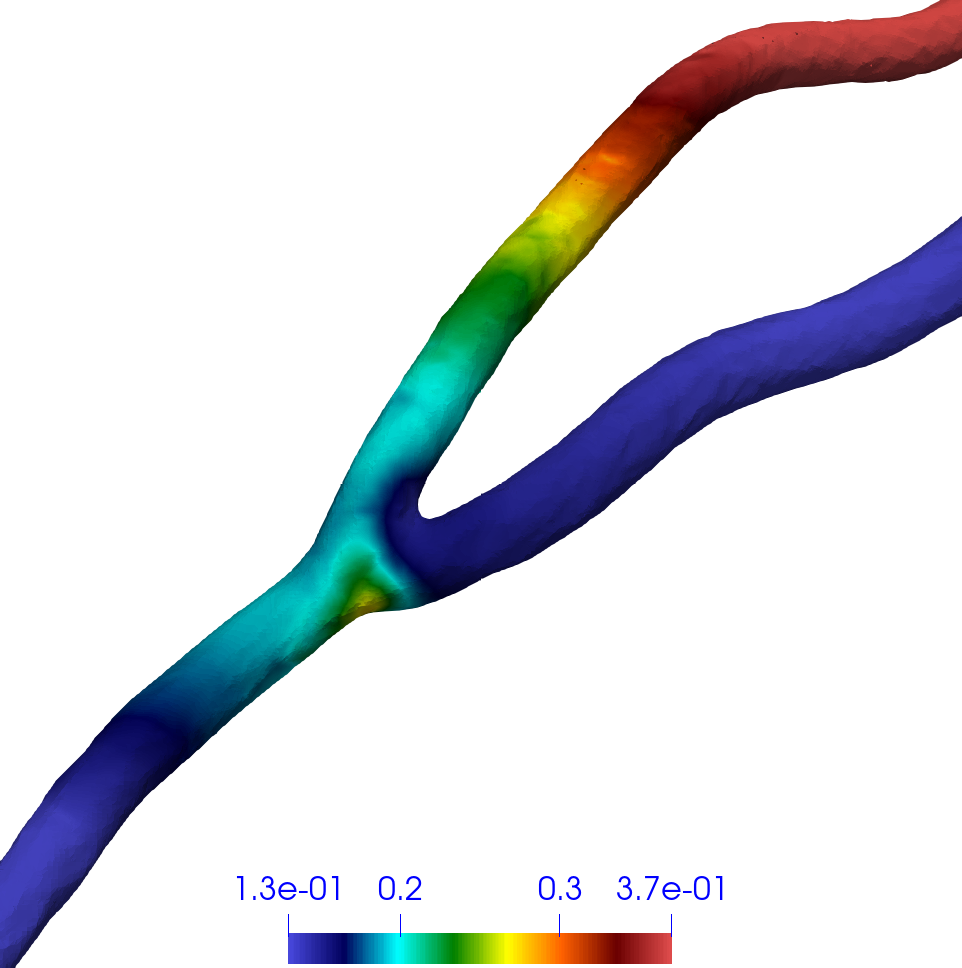}
			\put(20,105){ROM}
	\end{overpic}}

	\caption{Comparison between normalized pressure drop $P^* = P/P_{max}$ 
		in the anastomosis region computed by the FOM and by the ROM for $t/T = 0.8$. We set $\delta=0.99$.}
	\label{p_graft}
\end{figure}

\begin{figure}
	\centering
	\subfloat[][]{%
		\begin{overpic}[width=0.32\textwidth]{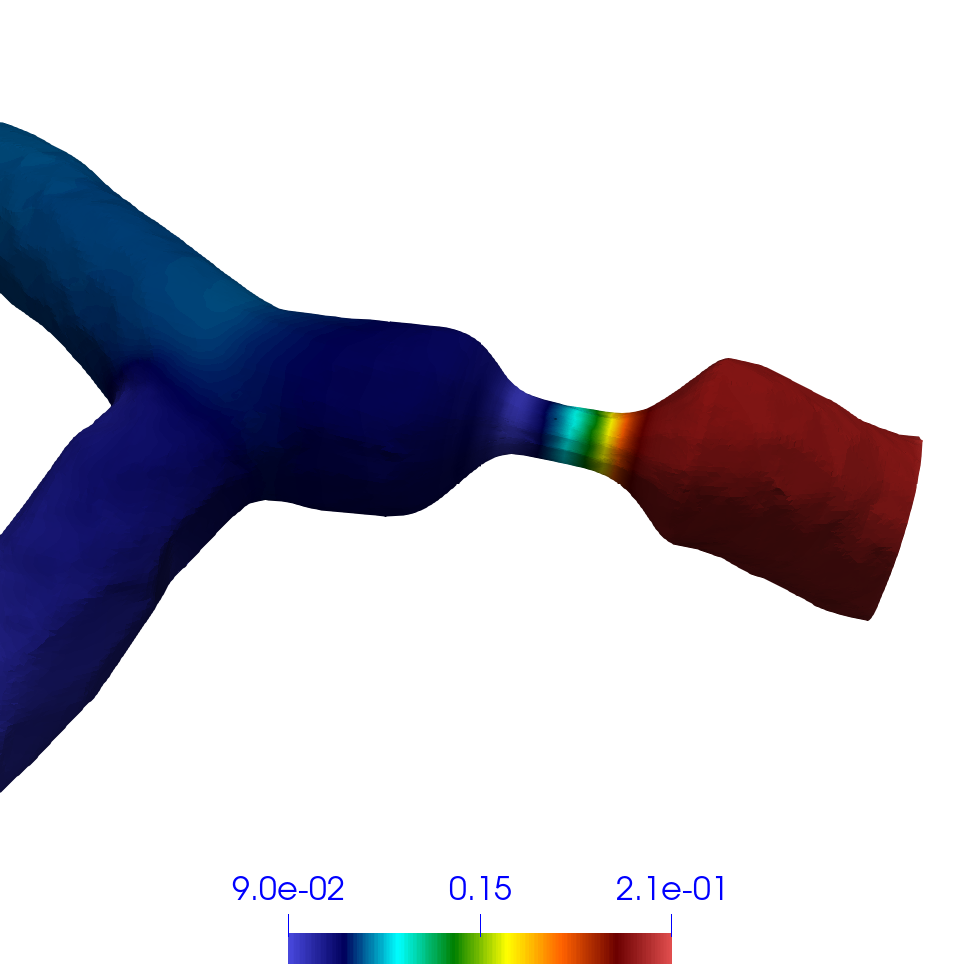}
			\put(20,105){FOM}
	\end{overpic}}
	\subfloat[][]{%
		\begin{overpic}[width=0.32\textwidth]{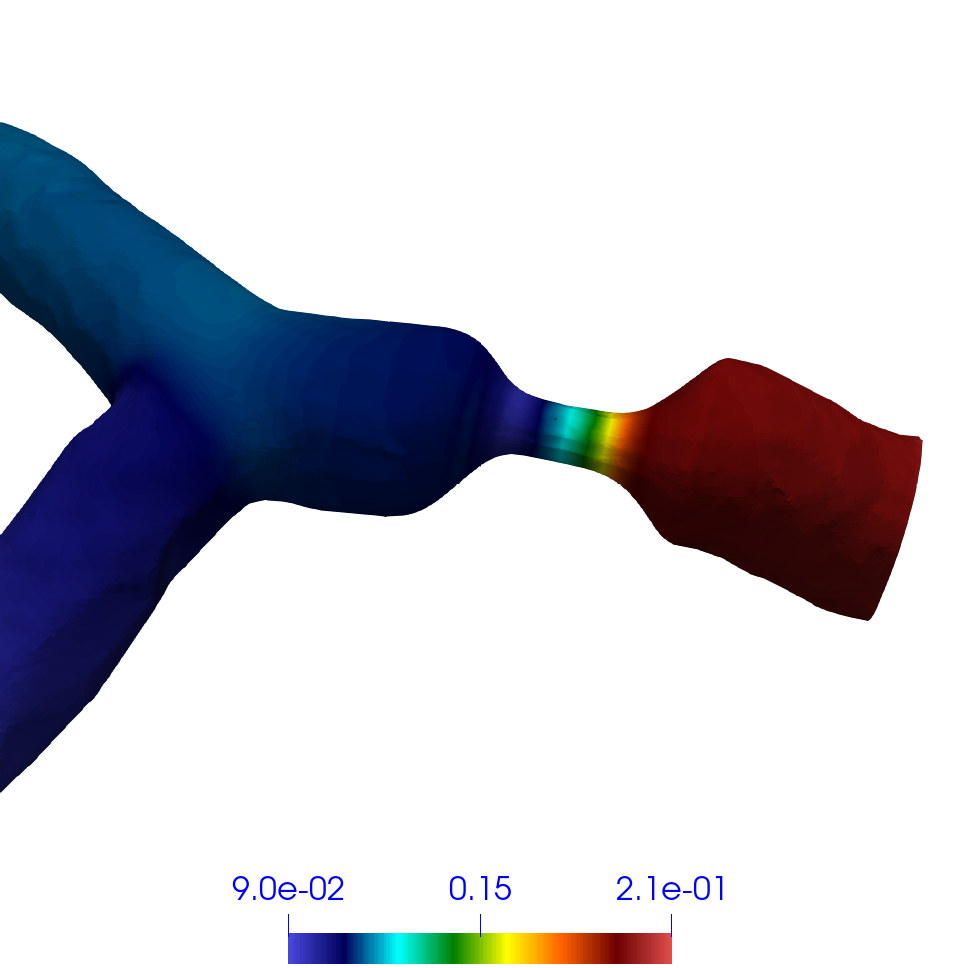}
			\put(20,105){ROM}
	\end{overpic}}
	\caption{Comparison between normalized pressure drop $P^* = P/P_{max}$ 
		in the stenosis region computed by the FOM and by the ROM for $t/T = 0.8$. We set $\delta=0.99$.}
	\label{p_stenosi}
\end{figure}

\begin{figure}
	\centering
	
	\subfloat[][]{%
	\begin{overpic}[width=0.32\textwidth]{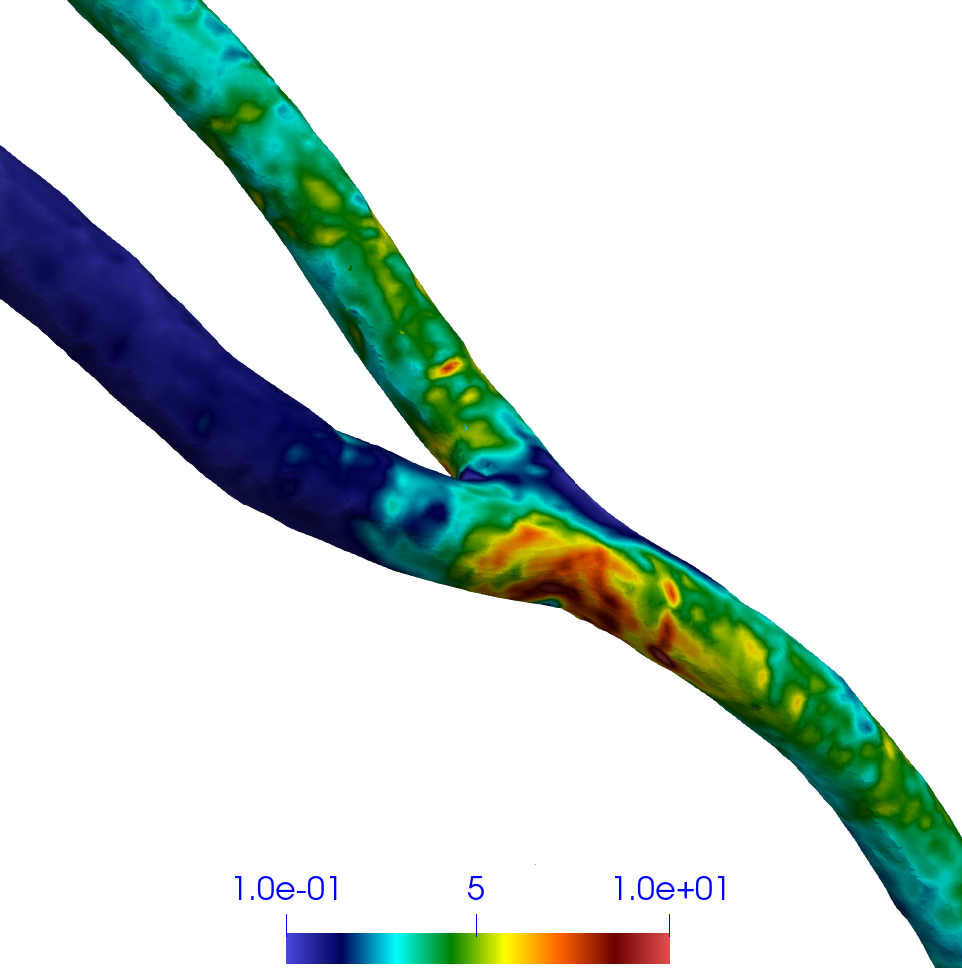}
			\put(20,105){FOM}
	\end{overpic}}
	\subfloat[][]{%
	\begin{overpic}[width=0.32\textwidth]{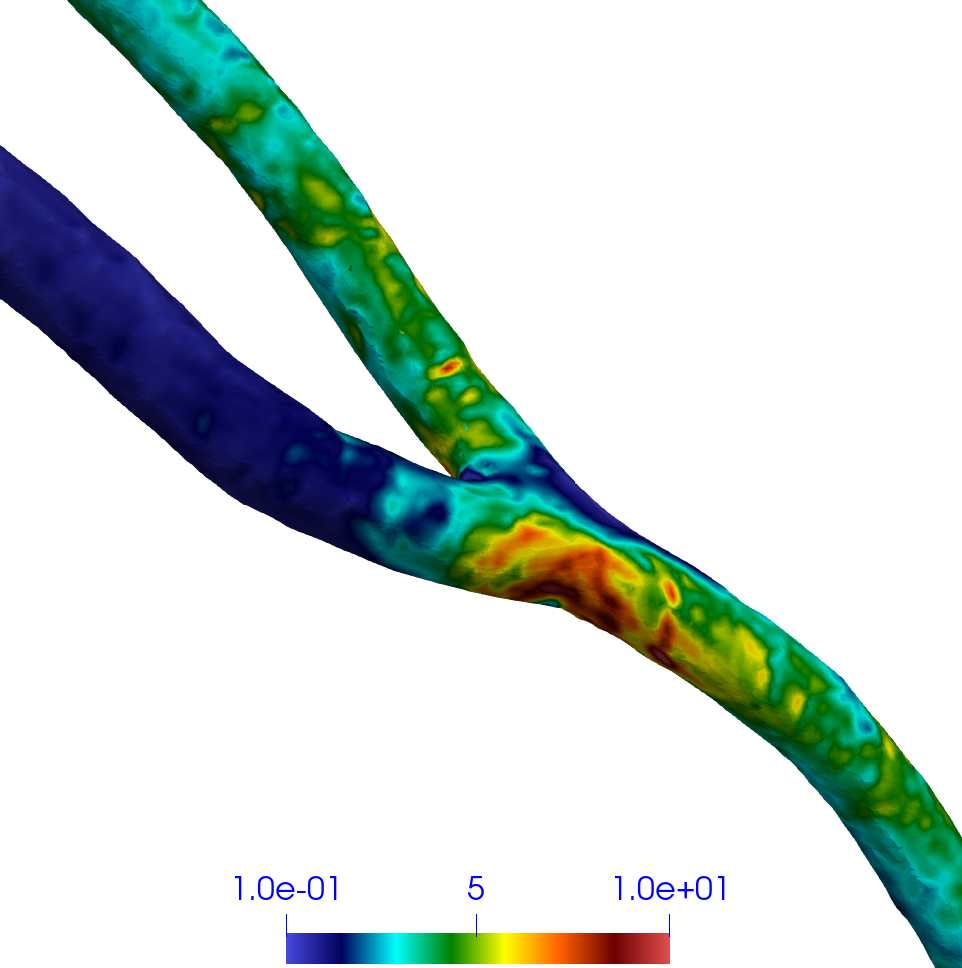}
			\put(20,105){ROM}
	\end{overpic}}
	\caption{Comparison between WSS distribution (Pa) in the anastomosis region computed by the FOM and by the ROM for $t/T = 0.8$. We set $\delta=0.99$.}
	\label{wss_graft}
\end{figure}

\begin{figure}
	\centering
	\subfloat[][]{%
	\begin{overpic}[width=0.32\textwidth]{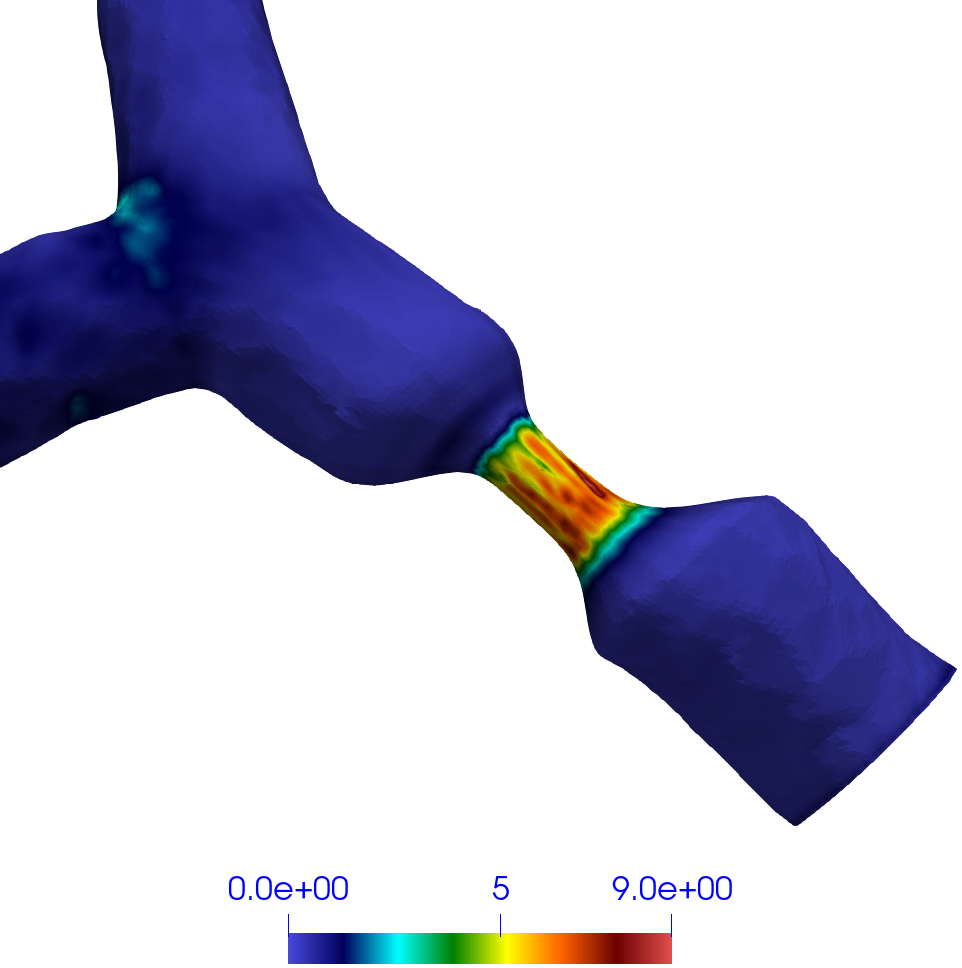}
			\put(20,105){FOM}
	\end{overpic}}
	\subfloat[][]{%
	\begin{overpic}[width=0.32\textwidth]{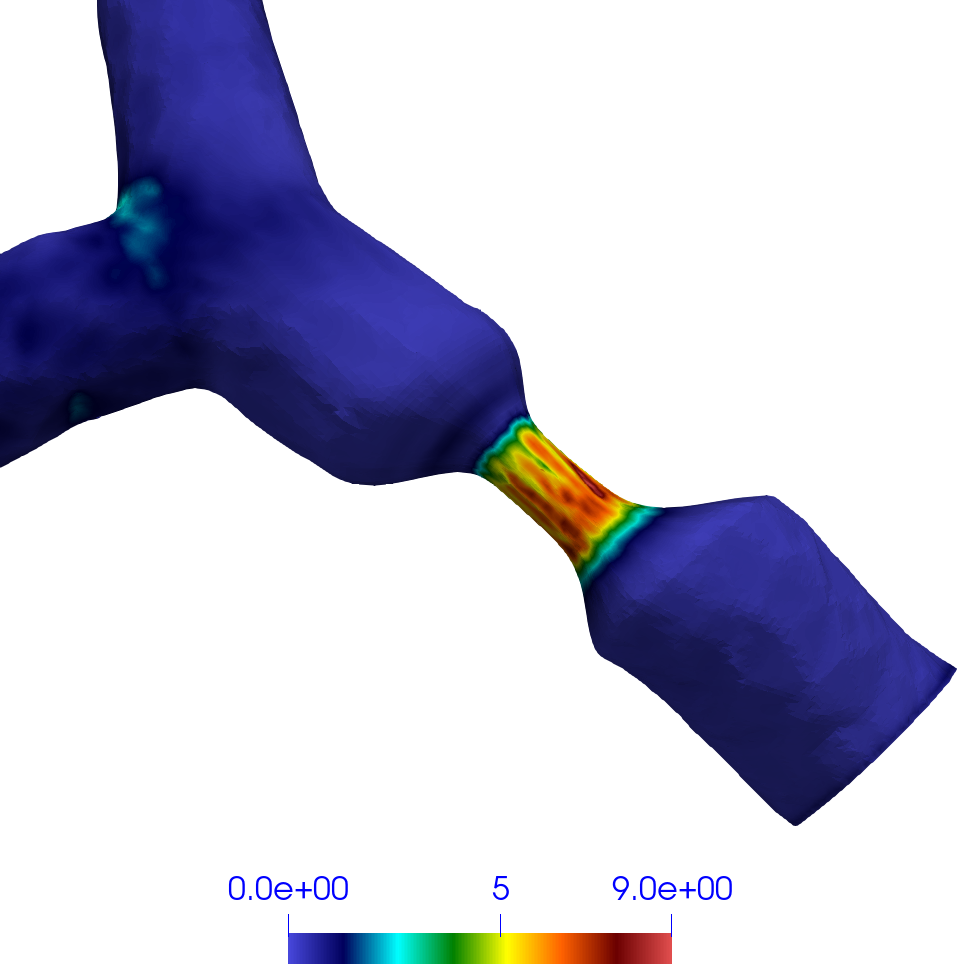}
			\put(20,105){ROM}
	\end{overpic}}
	\caption{Comparison between WSS distribution (Pa) in the stenosis region computed by the FOM and by the ROM for $t/T = 0.8$. We set $\delta=0.99$.}
	\label{wss_stenosi}
\end{figure}

\begin{figure}
	\centering
	\subfloat[][]{%
	\begin{overpic}[width=0.32\textwidth]{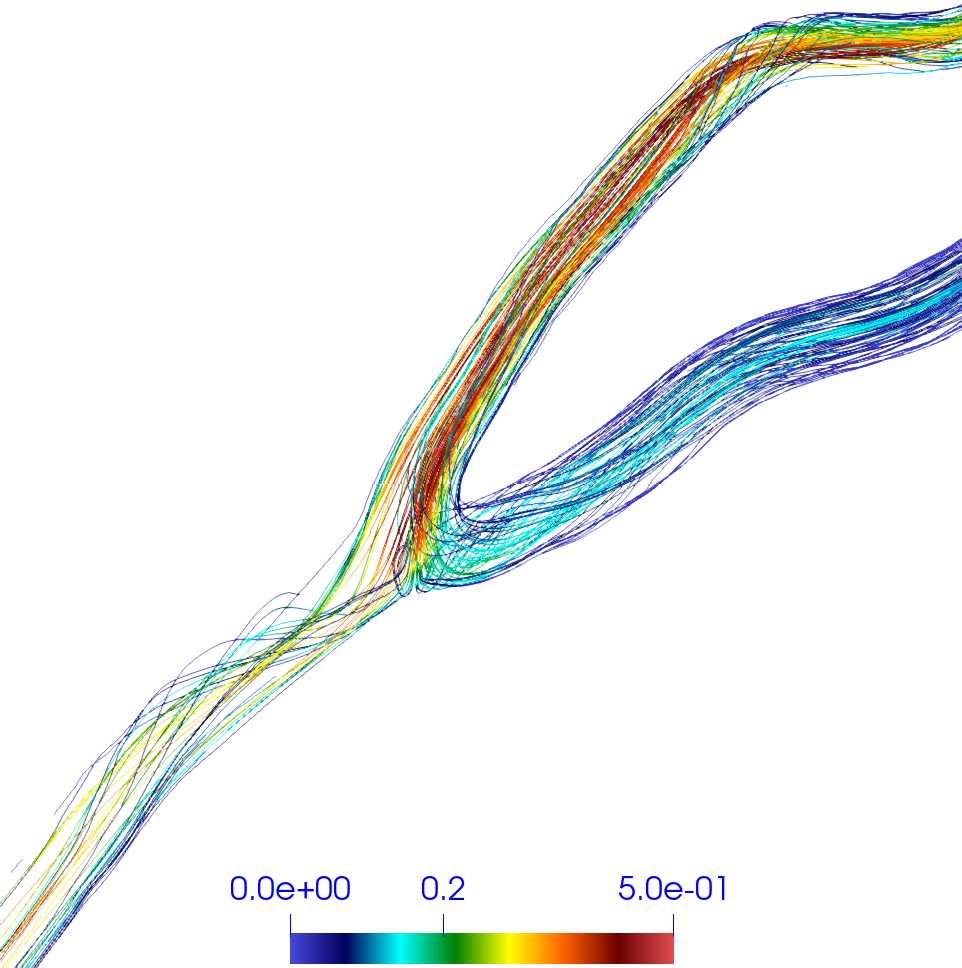}
			\put(20,105){FOM}
	\end{overpic}}
	\subfloat[][]{%
	\begin{overpic}[width=0.32\textwidth]{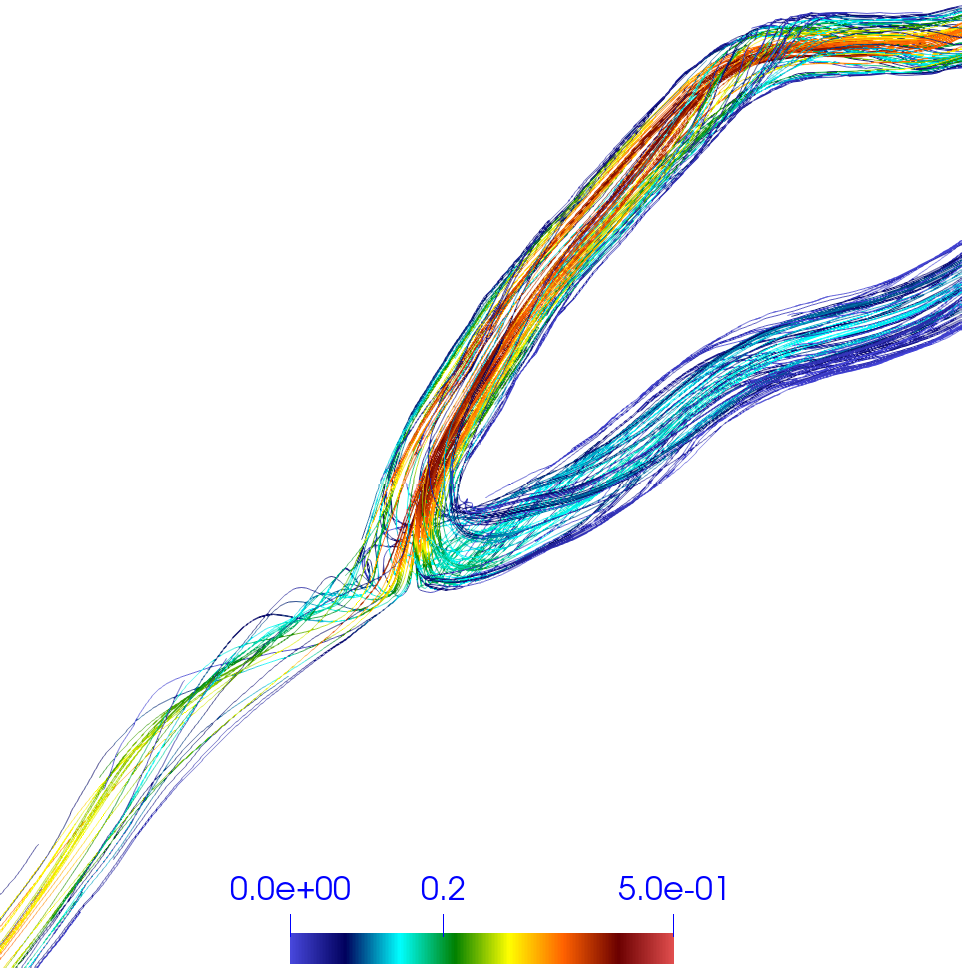}
			\put(20,105){ROM}
	\end{overpic}}
	\caption{Comparison between velocity streamlines (m/s) in the anastomosis region computed by the FOM and by the ROM for $t/T = 0.8$. We set $\delta=0.99$.}
	\label{stream_graft}
\end{figure}

\begin{figure}
	\centering
	\subfloat[][]{%
	\begin{overpic}[width=0.32\textwidth]{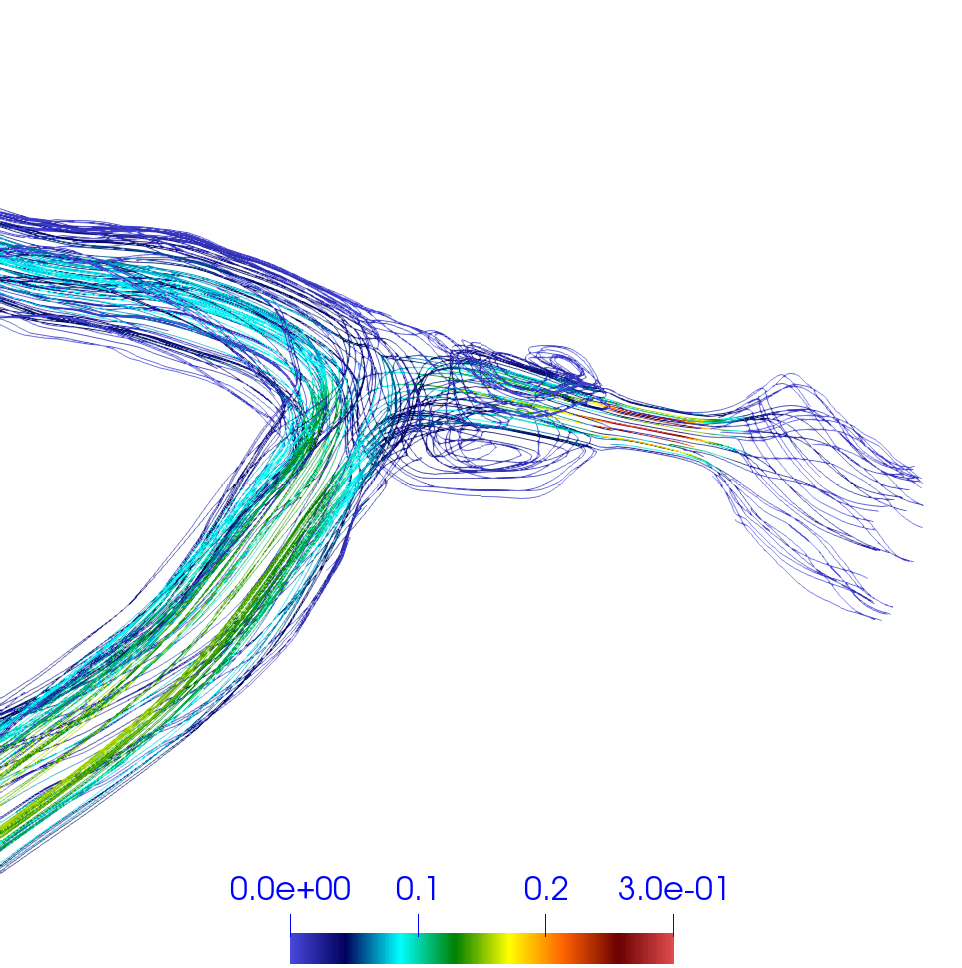}
		\put(20,105){FOM}
	\end{overpic}}
	\subfloat[][]{%
	\begin{overpic}[width=0.32\textwidth]{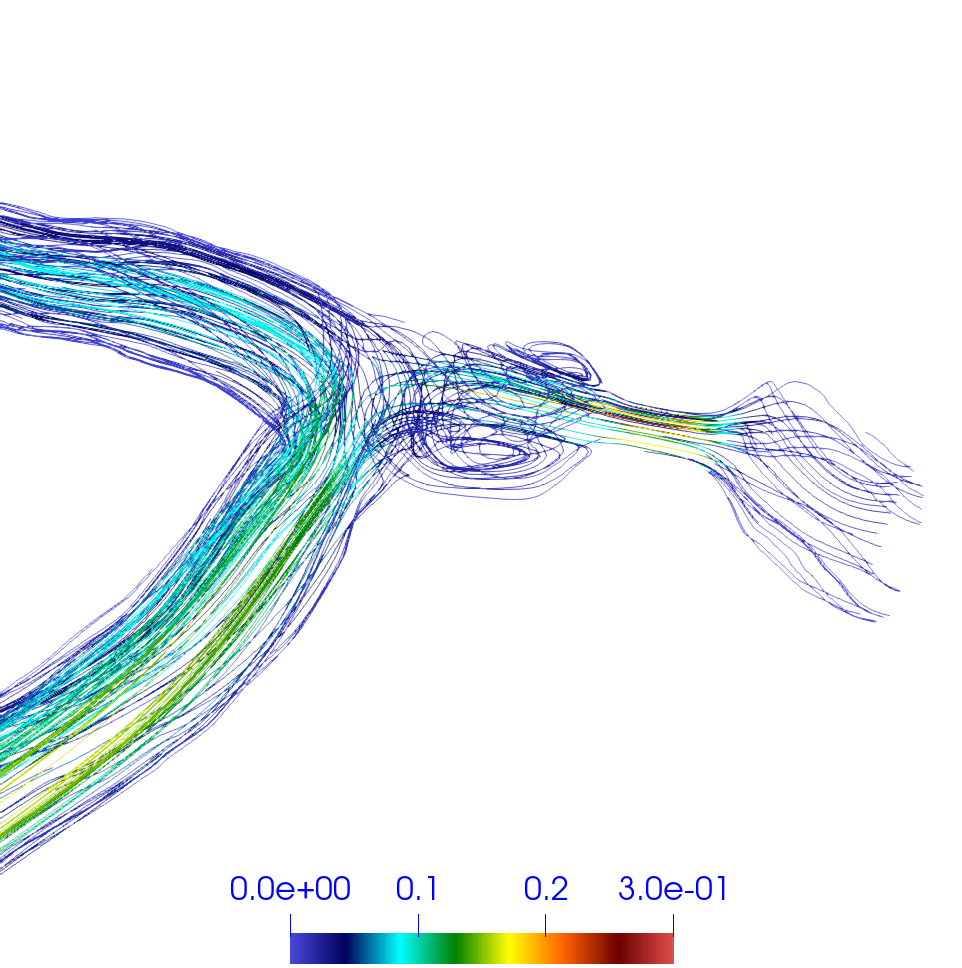}
		\put(20,105){ROM}
	\end{overpic}}
	\caption{Comparison between velocity streamlines (m/s) in the stenosis region computed by the FOM and by the ROM for $t/T = 0.8$. We set $\delta=0.99$.}
	\label{stream_stenosi}
\end{figure}

A qualitative comparison between FOM and ROM simulations is reported in {Figs.} \ref{p_graft} - \ref{stream_stenosi} for $t/T=0.8$. We display the stenosis and the anastomosis regions that, from the medical viewpoint, are those of major interest because they are modifications with respect to the healthy configuration. As one can see from {Figs.} \ref{p_graft} - \ref{stream_stenosi}, our ROM is able to provide a good reconstruction of velocity and pressure fields, as well as WSS field. 
{Figs.} \ref{p_graft} and \ref{p_stenosi} show the normalized pressure drop $P^* = P/P_{max}$ 
which is a useful indicator in the clinical practice to detect the presence of a stenosis and to measure its severity. 
{Figs.} \ref{wss_graft} and \ref{wss_stenosi} report the WSS distribution. As expected, a region of locally high WSS is found near the anastomosis and across the stenosis. It can represent a significant indication for the restenosis process.
Velocity streamlines are depicted in {Figs.} \ref{stream_graft} and \ref{stream_stenosi}. The velocity is higher in the LITA because it supplies blood to the entire vessels network, oxygenating LAD and, going up this vessel, LCx too. In {Fig.} \ref{stream_stenosi}, we observe that the velocity is elevated in the stenosis region due to the decreasing diameter. In addition, it can be seen a flow recirculation zone downstream of the stenosis. This phenomenon could favour the deposit of fatty materials on the surface of the lumen, causing the extension of the stenosis. 

Finally, we comment on the computational costs. We ran the FOM simulations in parallel using 20 processor cores. The simulations are run on the
SISSA HPC cluster Ulysses (200 TFLOPS, 2TB RAM, 7000 cores). Each FOM simulation takes roughly 41 h in terms of wall time, or 820 h in terms of total CPU time (i.e., wall time multiplied by the number of cores). On the other hand, the ROM has been run on an Intel(R) Core(TM) i5-8265U CPU @ 1.60GHz 8GB RAM by using one processor core only. Our ROM approach takes about 10 s for the computation of the reduced coefficients. So we obtain a speed-up of about $3 \cdot 10^{5}$. 


\section{PERSPECTIVES AND CONCLUSIONS}
\label{concluison}
In this work, a non-intrusive data-driven ROM based on POD-ANN approach is developed for fast and reliable numerical simulation of blood flow patterns occurring in a patient-specific
coronary system when an isolated stenosis of the LMCA occurs. A CABG performed with the LITA on the LAD is analyzed. The introduction of a patient-specific configuration is an attracting element of this work because it allow to establish a personalized clinical treatment. In addition, it is used a FFD technique, which gives the opportunity to deform directly the mesh and not only the geometry. Furthermore, the combination of ROM, FV technique and neural networks makes this study mathematically appealing. 

However, some insights are feasible. Parametric studies involving both physical and geometrical parameters will be performed in \citealp{pier2022neural}. 
Moreover, one could introduce Windkessel models (\citealp{girfoglio2020non,girfoglio2021non,fevola2021optimal}) in order to enforce more realistic outflow boundary conditions, which represent a crucial step to obtain meaningful outcomes. In addition, it could be interesting to investigate other ROM frameworks based on other deep learning approaches with the aim to furthermore improve both efficiency and accuracy of the method, such as physics-informed neural network (\citealp{chen2021physics,kissas2020machine,demo2021extended}).

Another important aspect deals with the combination of the ROM framework developed with technological progress through a web interface that would allow real time data to be accessed in hospitals and operating rooms on portable devices. In this scenario, the web server ARGOS (\url{https://argos.sissa.it}) has been created, which has the task of proposing a platform to favor a more widespread exploitation of real time
computing through a simple "click". ARGOS offers a wide variety of applications related to several problems, implemented by using boh intrusive (RBniCS, ITHACA-FV) and non intrusive (EZyRB, PyDMD) ROM-based softwares developed by SISSA mathLab (\url{https://mathlab.sissa.it/cse-software}), and in particular it contains the section ATLAS (\url{https://argos.sissa.it/atlas}, see \cite{girfoglio2021non} for further details) focused on the cardiovascular field. 


\section*{ACKNOWLEDGMENTS}
We acknowledge the support provided by the European Research Council Executive Agency by the Consolidator Grant project AROMA-CFD "Advanced Reduced Order Methods with Applications in Computational Fluid Dynamics" - GA 681447, H2020-ERC CoG 2015 AROMA-CFD, PI G. Rozza, and INdAM-GNCS 2019-2020 projects.


\bibliographystyle{elsarticle-harv}  
\bibliography{references.bib}


\end{document}